\documentclass[11pt,a4paper, parskip=half]{scrartcl} 
\usepackage[utf8]{inputenc}
\usepackage[english]{babel}
\usepackage[T1]{fontenc}
\usepackage{lmodern}
\usepackage[left=3cm,right=3cm,top=3cm,bottom=3cm]{geometry}
\usepackage[runin]{abstract}
    \abslabeldelim{.}

\usepackage{todonotes}

\usepackage{amsmath, amsfonts, amssymb, amsthm, amstext}
\usepackage{mathtools}
\usepackage{mathrsfs}
\usepackage{mathdots}
\usepackage{cases}

\usepackage{graphicx}
\usepackage{subcaption}
\usepackage{algorithm}

\font\mfett=cmmib10 at11pt
 at9pt
\def\bgamma{\hbox{\mfett\char013}}
\def\gamra{\hbox{\mfett\char013}}

\usepackage{enumitem}

\usepackage{hyperref} 
\usepackage[capitalise]{cleveref}
	
\newcounter{thm}
\numberwithin{thm}{section}
\numberwithin{equation}{section}

	\newtheoremstyle{myplain}		
			{}			
			{}			
			{\itshape}				
			{}				
			{\sffamily\bfseries}				
			{.}		
			{ }				
			{\thmname{#1}\thmnumber{ #2}\textnormal{\textsf{\thmnote{ (#3)}}}}			
    \newtheoremstyle{mybreak}
            {}{}{}{}{\sffamily\bfseries}{.}{\newline}
            {\thmname{#1}\thmnumber{ #2}\textnormal{\textsf{\thmnote{ (#3)}}}}
	\newtheoremstyle{mydef}
			{}{}{}{}{\sffamily\bfseries}{.}{ }
			{\thmname{#1}\thmnumber{ #2}}
	\newtheoremstyle{myrem}
			{}{}{}{}{\sffamily\itshape}{.}{ }
			{\thmname{#1}\thmnumber{ #2}}

\theoremstyle{myplain}

\theoremstyle{mybreak}
\theoremstyle{mydef}
	
	\newtheorem{remark}[thm]{Remark}
\theoremstyle{mydef}
	\newtheorem{example}[thm]{Example}

	\newcommand{\cc}{\mathbb{C}}
		\newcommand{\zz}{\mathbb{Z}}
		\newcommand{\nn}{\mathbb{N}}
	\newcommand{\rr}{\mathbb{R}}

\DeclareMathOperator{\diag}{diag}

\allowdisplaybreaks

\def\sumprime_#1^#2{
    \setbox0=\hbox{$\scriptstyle{#1}$}
    \setbox1=\hbox{$\scriptstyle{#2}$}
    \setbox2=\hbox{$\displaystyle{\sum}$}
    \setbox4=\hbox{${}^\prime\mathsurround=0pt$}
    \dimen0=.5\wd0 \advance\dimen0 by-.5\wd2
    \ifdim\dimen0>0pt
        \ifdim\dimen0>\wd4 \kern\wd4
        \else\kern\dimen0
        \ifdim\dimen1>\wd4 \kern\wd4
        \else\kern\dimen1
    \fi\fi\fi
\mathop{{\sum}^\prime}_{\kern-\wd4 #1}^{\kern-\wd4 #2}
}

\title{\Large Recovery of rational  functions  via Hankel  pencil method and sensitivities of the poles}
\author{Nadiia Derevianko\footnote{TUM School of CIT, 
Department of Computer Science,
Boltzmannstrasse 3,
85748 Garching b. München,
Germany,  nadiia.derevianko@tum.de} \footnote{Corresponding author} }
\date{\today}

\begin{document}
	\let\oldproofname=\proofname
	\renewcommand{\proofname}{\itshape\sffamily{\oldproofname}}

\maketitle

\begin{abstract}
In this paper, we introduce a new approach for the recovery of rational functions. The concept we propose is based on using  the exponential structure of the Fourier coefficients of rational functions and the reconstruction of this exponential structure in the frequency domain. We choose ESPRIT as a method for the exponential recovery. The matrix pencil structure of this approach is the reason for its selection, as it makes our method suitable for the sensitivity analysis.  According to our method, poles located inside and outside the unit circle are reconstructed independently as eigenvalues of some special Hankel matrix pencils.  
Furthermore, we derived formulas for sensitivities of poles of rational functions  in case of unstructured and structured perturbations. Finally, we consider several numerical experiments and,  using sensitivities,  explain the recovery errors for poles.
\\[1ex]
\end{abstract}

\section{Introduction }

For $M \in \nn$,  we define \textit{a rational function}
 $r(z)$  of type $(M-1,M)$ by means of its partial fraction decomposition
\begin{equation}\label{rat}
r(z):=\sum\limits_{j=1}^{M}\frac{\gamma_j}{z-z_j},
\end{equation}
where $\gamma_j,  \,  z_j \in \cc$, coefficients $\gamma_j  \neq 0$ and poles $z_j$ are pairwise distinct.    
 We study the question of recovery of rational functions (\ref{rat}) using as input data a finite set of their Fourier coefficients $\widehat{r}_k$ defined by (\ref{rkf}).   We use the property that Fourier coefficients of rational functions have the exponential structure, and we reconstruct this exponential structure in the frequency domain applying  the ESPRIT method \cite{RK89,  PT2013}.  The idea of recovery of rational functions by estimation of signal parameters in the frequency domain is not new.  
Similar technique for the recovery of or approximation by rational functions  was employed in  \cite{DB13,  WDT21, Y22,  Y222}. The paper \cite{Y22} is of special interest to us. In \cite{Y22}, the author considered recovery of  rational functions (\ref{rat}) based on using noisy input data  and discussed reconstruction results for the poles $z_j$, $j=1,...,M$. It was pointed out that for different levels of noise, not all poles $z_j$ are reconstructed with the same accuracy, and sometimes even the locations of certain poles are not reconstructed  correctly. From these numerical observations, we can conclude  that some poles $z_j$ of  $r$ as in (\ref{rat}) are more ''sensitive`` to perturbations in the input data than the others.  Our goal in this paper is to present theoretical explanations of these numerical observations by applying the \textit{sensitivity analysis}. The method developed in \cite{Y22} employees the  AAK theory \cite{PP}  for the reconstruction of the exponential structure of the Fourier coefficients $\widehat{r}_k$. We modify these ideas,  by applying the ESPRIT method in the frequency domain.  Since ESPRIT is based on using the Hankel pencil,  it makes our method suitable for the sensitivity analysis.

In \cite{ZGA21}, the authors developed the sensitivity analysis for the eigenvalues  of Loewner pencil  with respect to the unstructured and structured perturbations. We refer the interested reader to \cite{ALI17}  for further application of Loewner framework for model reduction of dynamical systems. 
 Following the approach in \cite{ZGA21},  we present the sensitivity analysis for the poles of $r$ in (\ref{rat}) computed by our method.  According to our idea,  poles located inside  and outside the unit circle are reconstructed separately as eigenvalues of some Hankel matrix pencils.  Further,  we consider two kinds of Hankel pencil perturbations.  For the sensitivity $\rho$ with respect to \textit{unstructured perturbations},  we present the connection with the condition numbers of Vandermonde matrices that are used in factorizations of the corresponding Hankel matrices.  The sensitivity $\eta$ with respect to \textit{structured perturbations} of Hankel pencil is studied in \cite{Zhphd}.  We present modification of these results for our special rational reconstruction problem (\ref{rat}).   Note that the theory of sensitivities of eigenvalues  of Hankel pencil perturbations was considered also in \cite{G99, BGL}.  In \cite{G99} and  \cite{BGL},  the authors measured both structured and unstructured perturbations by the same quantity.  Similarly as in  \cite{ZGA21},   we introduce two different quantities, sensitivities $\eta$ and  $\rho$ with respect to structured und unstructured perturbations, respectively, for poles of the rational function   (\ref{rat}).  The ill-disposedness of $z_j$ is measured by determining whether $\rho$ and $\eta$ are ''large`` for this specific pole $z_j$. Additionally, the sensitivity $\eta$ reveals the perturbations of which Fourier coefficient $\widehat{r}_k$ the pole $z_j$ is most sensitive to.

The method developed in this  paper is also closely related to our ESPIRA method for exponential recovery \cite{DPP21}--\cite{DP21} (see also \cite{DH25} for the multivariate case).  ESPIRA uses the property that the classical Fourier coefficients of exponential sums have a rational structure and reconstructs this rational structure with the AAA method for rational approximation \cite{AAA}.  Thus, the rational reconstruction concept employed in this paper can be viewed as the opposite of the ESPIRA concept. We would also like to draw particular attention to the paper \cite{BC20}, in which the authors used the sensitivity analysis of Hankel pencil eigenvalues to developed a method for exponential reconstruction.

Finally, we want to mention the other methods for the recovery of rational functions such as the AAA method \cite{AAA},   the RKFIT (Rational Krylov Fitting) method  \cite{BG17},  the Vector Fitting  algorithm  \cite{GS99} and the  Loewner framework method \cite{ ALI17}. In  \cite{GG21}, the authors discussed all these algorithms for the recovery of rational functions using as input information given noisy set of samples of $r$.  Since we are interested in the sensitivity analysis of the poles, we have to assume that we apply perturbations to the Fourier coefficients $\widehat{r}_k$. Note that in this paper, we also introduce the corresponding modification of our method, which is based on  employing function values of $r$ by computing $\widehat{r}_k$ via a quadrature formula for some grid of points (see Remark \ref{rem41}).  Modifying our algorithm in this way,  allows us to compare its performance with known rational methods. Numerical experiments show us that our approach gives equally good and, for some examples, even better recovery results from given noisy function values than the AAA method.

\textbf{Outline. }  In Section 2,  we present the main ideas of the ESPRIT method for the exponential recovery.   In Section 3, we consider the sensitivity analysis for Hankel pencil eigenvalues:  Subsection 3.1 is dedicated to unstructured perturbations and in Subsection 3.2,  we consider structured perturbations.  In Section 4,  we describe our main results.  In Subsection 4.1,  we present main ideas of our new approach for the recovery of rational functions (\ref{rat}).  Subsection 4.2 is dedicated to the sensitivity analysis of poles reconstructed with our method.   Finally,  in Section 5,  we present numerical experiments.

\textbf{Notation.} As usual, $\nn$, $\zz$, $\rr$ and $\cc$ denote the natural, integer,  real and complex numbers, respectively. 
Let us agree that throughout the paper,  we  use the matrix notation $\mathbf{X}_{L,K}:=(X_{i,j})_{i,j=1}^{L,K} $ for   matrices of the size $L\times K$ and the submatrix notation $\mathbf{X}_{L,K}(m:n,k:\ell)$,  to denote  a submatrix of $\mathbf{X}_{L,K}$ with rows indexed from $m$ to $n$ and columns indexed from $k$ to $\ell$,  where first row and first column has index 1 (even though the row and column indices for the definition of the matrix may start with 0).  For square matrices,  we often use the short notation $\mathbf{X}_{K}$ instead of $\mathbf{X}_{K,K}$.  By $ \mathbf{X}_{L,K}^{\dagger}$ we denote the Moore-Penrose inverse of  $\mathbf{X}_{L,K}$ and by   $ \mathbf{X}_{K}^{-1}$ inverse of the invertible matrix $\mathbf{X}_{K}$.   We consider $\|\mathbf{X}_{L,K}\|_2:=\sigma_{\max}(\mathbf{X}_{L,K})$ the 2-norm or spectral norm,  and $\|\mathbf{X}_{L,K}\|_{\mathrm{F}}:= \left(\sum_{i=1}^{L} \sum_{j=1}^{K}   |X_{i,j}|^{2} \right)^{1/2}=\left(\sum_{i=1}^{\min\{L,K \}} \sigma_i^{2}  \right)^{1/2}$ the Frobenius norm of the matrix $\mathbf{X}_{L,K}$. By $\boldsymbol{x}=(x_1,\ldots,x_K)^{T}$ we denote a column vector of the size $K$.  The $\ell_1$- and $\ell_2$-norms of a vector $\boldsymbol{x}$ are given  by $\|\boldsymbol{x}\|_1=\sum_{j=1}^{K} |x_j|$ and $\|\boldsymbol{x}\|_2=\left( \sum_{j=1}^{K} |x_j|^{2} \right)^{1/2}$, respectively.

\section{ESPRIT method for the recovery of exponential sums}
\label{esprit}

As previously mentioned in Introduction, we employ the ESPRIT method \cite{RK89, PT2013} to reconstruct the exponential structure of the Fourier coefficients of rational functions (\ref{rat}). 
In this section,  we describe the key concepts of this approach.  

Let  $f$ be \textit{an exponential sum} of order $M\in \nn$ defined by
\begin{equation}\label{expsum}
f(t):=\sum\limits_{j=1}^{M} \gamma_j z_j^{t},   \  \gamma_j,  z_j \in \cc, 
\end{equation}
where nodes $z_j$ are pairwise distinct and coefficients $\gamma_j \neq 0$.
Let $L, N \in \nn$ be given such that  $M\leq L  \leq N$,  and $L$ is the upper bound for $M$.  Our goal is to reconstruct the exponential sum (\ref{expsum}),  i.e.  to reconstruct parameters $M$,  $\gamma_j$ and $z_j$ for $j=1,\ldots,M$,  using as input data a vector of  samples $\boldsymbol{f}:  = \big(f(0),  f(1), \ldots, f(2N-1)  \big)^{T}$ of the size $2N$.  

 For $n,m \in \nn$ and  a vector $\boldsymbol{a}:=(a_1,\ldots, a_m)^{T}$,  we define a Vandermonde  matrix $\mathbf{V}_{n,m}: =\mathbf{V}_{n,m}(\boldsymbol{a})=(a_j^{k})_{k=0,j=1}^{n-1,m}$.  Let now  $\boldsymbol{z}:=(z_1,\ldots, z_M)^{T}$ be a vector of nodes of the exponential sum (\ref{expsum}) and let  $\mathbf{\Gamma}_M:=\diag(\gamma_j)_{j=1}^{M}$ and $\mathbf{Z}_M:=\diag(z_j)_{j=1}^{M}$.  We create a Hankel matrix ${\mathbf H}_{2N-L, L+1} \coloneqq \left( f(\ell+m)\right)_{\ell,m=0}^{2N-L-1,L}$ and consider its factorization 
\begin{equation}\label{factH}
\mathbf{H}_{2N-L,L+1}=\mathbf{V}_{2N-L,M} (\boldsymbol{z}) \,  \mathbf{\Gamma}_M \,  \mathbf{V}^{T}_{L+1,M} (\boldsymbol{z}).
\end{equation}
   Since Vandermonde  matrices in (\ref{factH}) have full column ranks,  the matrix $\mathbf{H}_{2N-L,L+1}$ has rank $M$.  Thus, in case of exact given data,  the order $M$ of the exponential sum (\ref{expsum})   is reconstructed as the numerical rank of the matrix $   \mathbf{H}_{2N-L,L+1}$.
 
 Consider further a Hankel matrix pencil problem   
 \begin{equation}\label{1.6}
z \mathbf{H}_{2N-L,L}(0) -  \mathbf{H}_{2N-L,L}(1),
\end{equation}
where  $\mathbf{H}_{2N-L,L}(s):=\mathbf{H}_{2N-L,L+1}(1:2N-L,1+s:L+s)=(f(\ell+k+s))_{k,\ell=0}^{2N-L-1,L-1}$ for $s=0,1$.
   Since the  Hankel matrices $\mathbf{H}_{2N-L,L}(0)$ and $\mathbf{H}_{2N-L,L}(1)$ satisfy factorizations
\begin{align}
\mathbf{H}_{2N-L,L}(0) &=\mathbf{V}_{2N-L,M} (\boldsymbol{z}) \,  \mathbf{\Gamma}_M \,  \mathbf{V}^{T}_{L,M} (\boldsymbol{z}), \label{fachan} \\
\mathbf{H}_{2N-L,L}(1)&=\mathbf{V}_{2N-L,M} (\boldsymbol{z}) \,  \mathbf{Z}_M \,  \mathbf{\Gamma}_M \,  \mathbf{V}^{T}_{L,M}(\boldsymbol{z}), \label{fachan1}
\end{align}
we obtain that
$$ z \,  {\mathbf H}_{2N-L,L}(0) - {\mathbf H}_{2N-L,L}(1) = \mathbf{V}_{2N-L,M} (\boldsymbol{z})   \, \mathrm{diag}   \left( z - z_{j}\right)_{j=1}^{M}  \,  \mathbf{\Gamma}_M \,  \mathbf{V}^{T}_{L,M}(\boldsymbol{z}).$$
 From the last formula,  we conclude that nodes $z_j$,  $j=1,\ldots,M$   in (\ref{expsum})  are eigenvalues of the matrix pencil (\ref{1.6}).    In case when $M$ is known,  we can choose $N=L=M$ and compute $z_j$ as eigenvalues of 
  \begin{equation}\label{mp11}
z \mathbf{H}_{M}(0) -  \mathbf{H}_{M}(1).
\end{equation}

Next, we discuss how we solve the matrix pencil (\ref{1.6}).
  Consider  SVD of  $ {\mathbf H}_{2N-L,L+1}$,
\begin{equation}\label{SVD}
 {\mathbf H}_{2N-L,L+1} = {\mathbf U}_{2N-L} \, {\mathbf D}_{2N-L,L+1} {\mathbf W}_{L+1}, 
\end{equation} 
where  ${\mathbf U}_{2N-L}$ and ${\mathbf W}_{L+1}$  are unitary  matrices  and, in  case of unperturbed data,
$$ {\mathbf D}_{2N-L,L+1} = \left( \begin{array}{cc} \mathrm{diag} \left(  (\sigma_{j})_{j=1}^{M} \right) & {\mathbf 0} \\
{\mathbf 0} & {\mathbf 0} \end{array} \right) \in {\mathbb R}^{(2N-L) \times (L+1)}, $$
where $\sigma_{1} \ge  \ldots  \ge \sigma_{M}>0$ are singular values of ${\mathbf H}_{2N-L,L+1}$.
Consequently, (\ref{SVD}) can be simplified to 
\begin{equation}\label{1.7}
{\mathbf H}_{2N-L,L+1} = {\mathbf U}_{2N-L} \, {\mathbf D}_{2N-L,M} {\mathbf W}_{M,L+1}, 
\end{equation}
where the last $L+1-M$ zero-columns in ${\mathbf D}_{2N-L,L+1}$ and the last $L+1-M$ zero-rows in ${\mathbf W}_{L+1}$ are removed.
Similarly as before,  SVDs of ${\mathbf H}_{2N-L,L}(0)$ and ${\mathbf H}_{2N-L,L}(1)$ can be directly  derived  from (\ref{1.7})  and the matrix pencil (\ref{1.6}) now reads
\begin{eqnarray*}
 & & z \left(  {\mathbf U}_{2N-L} \, {\mathbf D}_{2N-L,M} {\mathbf W}_{M,L+1} \right)(1:2N-L,1:L)  \\
 & &  \hspace{5cm} -  \left(  {\mathbf U}_{2N-L} \, {\mathbf D}_{2N-L,M} {\mathbf W}_{M,L+1} \right)(1:2N-L,2:L+1),
\end{eqnarray*}
or equivalently
$$ z \left(   {\mathbf D}_{2N-L,M} {\mathbf W}_{M,L+1} \right)(1:2N-L,1:L) -  \left(  {\mathbf D}_{2N-L,M} {\mathbf W}_{M,L+1} \right)(1:2N-L,2:L+1). $$
Multiplication from the left with $\left( \mathrm{diag}( \sigma_{j}^{-1} )_{j=1}^{M} , {\mathbf 0} \right) \in {\mathbb R}^{M \times (2N-L)}$  yields
the matrix pencil
\begin{equation}\label{mp}
z \,  {\mathbf W}_{M,L+1}(1:M,1:L) - {\mathbf W}_{M,L+1}(1:M,2:L+1). 
\end{equation}
Therefore,  eigenvalues of the matrix pencil (\ref{1.6}) can be computed as eigenvalues of the matrix pencil (\ref{mp}).  Then we define ${\mathbf W}_{M,L}(0) \coloneqq {\mathbf W}_{L+1}(1\!:\!M,1\!:\!L)$ and
  ${\mathbf W}_{M,L}(1) \coloneqq {\mathbf W}_{L+1}(1\!:\!M,2\!:\!L+1)$,
and determine the vector $\boldsymbol{z} =(z_{1}, \ldots , z_{M})^{T}$ as a vector of eigenvalues of 
$\left(  {\mathbf W}_{M,L}(0)^{T} \right)^{\dagger}  {\mathbf W}_{M,L}(1)^{T}$.

Having nodes $z_j$,  we compute coefficients $\gamma_j$,  $j=1,\ldots,M$ as a least squares solution of the system of linear equations
$$
\sum\limits_{j=1}^{M} \gamma_j z_j^{k}=f(k),  \quad  k=0,\ldots,2N-1,
$$
which can be also written in a matrix form
$$
{\mathbf V}_{2N,M} ({\boldsymbol z}) \, \bgamma = \boldsymbol{f},
$$
where ${\bgamma}: = (\gamma_1,\ldots, \gamma_M)^{T}$.
The arithmetical complexity of the ESPRIT method is governed by the computation of SVD of $H_{2N-L,L+1}$ with $\mathcal{O}(N L^{2})$ operations.

\section{Sensitivities of Hankel pencil  eigenvalues}

 The concept of sensitivities arises from the perturbations of eigenvalues of the matrix pencil and    for Hankel matrices  was introduced in \cite{G99} and further considered in \cite{BGL}.   In  \cite{ZGA21},  the notion of sensitivities was defined for Loewner matrix pencil.  In the paper \cite{ZGA21},  the authors define two sensitivities,  which correspond to structured and unstructured perturbations.  Using similar ideas,  in \cite{Zhphd},  the author considered the structured  sensitivity of Hankel matrix pencil. In  \cite{BC20}, the   sensitivity analysis of Hankel pencil was applied to develop a new method for the exponential analysis. In this section,  we consider sensitivities of eigenvalues of the Hankel pencil (\ref{1.6}) and in Subsection \ref{sub42}, we present modification of these results  to develop sensitivity analysis for the poles of rational functions (\ref{rat}).

\subsection{Sensitivities of Hankel  pencil eigenvalues with unstructured perturbations}
\label{unsen}

Following the approach  for Loewner pencil from   \cite{ZGA21}, 
we show that in case of unstructured perturbations of Hankel pencil,  the sensitivities of eigenvalues depend on the numerical conditions of Vandermonde  matrices from factorizations (\ref{fachan})  and (\ref{fachan1}) (for $N=L=M$).
We consider the generalized eigenvalue problem (GEP) for the Hankel matrix pencil (\ref{mp11})
\begin{equation}\label{gep}
\boldsymbol{p}^{T} \mathbf{H}_{M}(1) = z \boldsymbol{p}^{T}\mathbf{H}_{M}(0),  \quad \mathbf{H}_{M}(1) \boldsymbol{q}= z \mathbf{H}_{M}(0) \boldsymbol{q} ,
\end{equation}
where $\boldsymbol{p}=(p_{1},\ldots,p_{M})^{T}$ and $\boldsymbol{q}=(q_{1},\ldots,q_{M})^{T}$ are left and right eigenvectors, respectively, of the GEP (\ref{gep}).
Then under the perturbations
\begin{equation}\label{pertsys}
\tilde{ \mathbf{H}}_{M}(0):=\mathbf{H}_{M}(0)+\Delta_{\mathbf{H}_{M}(0)},  \ \ \tilde{\mathbf{H}}_{M}(1):=\mathbf{H}_{M}(1)+\Delta_{\mathbf{H}_{M}(1)},
\end{equation}
 the first order approximation of the eigenvalue perturbation $z^{(1)}$ is given by 
\begin{equation}\label{pert}
z^{(1)}=\frac{\boldsymbol{p}^{T}\,  (\Delta_{\mathbf{H}_{M}(1)}-z \, \Delta_{\mathbf{H}_{M}(0)}) \,  \boldsymbol{q}}{\boldsymbol{p}^{T}  \,  \mathbf{H}_{M}(0) \,  \boldsymbol{q}}.
\end{equation}
Formula (\ref{pert}) is derived from the perturbed system 
\begin{align*}
(\boldsymbol{p}^{T}+(\tilde{\boldsymbol{p}}^{(1)})^{T}+\ldots) & (\mathbf{H}_{M}(1)+\Delta_{\mathbf{H}_{M}(1)}) = (z+z^{(1)}+\ldots) (\boldsymbol{p}^{T}+(\tilde{\boldsymbol{p}}^{(1)})^{T}+\ldots)(\mathbf{H}_{M}(0)+\Delta_{{\mathbf{H}_{M}(0)}}),\\
(\mathbf{H}_{M}(1)+\Delta_{\mathbf{H}_{M}(1)}) )& (\boldsymbol{q}+ \tilde{\boldsymbol{q}}^{(1)}+\ldots )= (z+z^{(1)}+\ldots)( \mathbf{H}_{M}(0)+\Delta_{\mathbf{H}_{M}(0)})  )(\boldsymbol{q}+ \tilde{\boldsymbol{q}}^{(1)}+\ldots ),
\end{align*}
by retaining only first-order terms (see \cite{G99},  \cite{BGL},  \cite{ZGA21} and \cite{Zhphd} for details).

Assume that the $\ell_2$-norms of perturbation matrices $\Delta_{\mathbf{H}_{M}(0)}$ and $\Delta_{\mathbf{H}_{M}(1)}$ in (\ref{pertsys}) are bounded as follows
\begin{equation}\label{normmat1}
\| \Delta_{\mathbf{H}_{M}(0)} \|_2\leq \epsilon \| \mathbf{H}_{M}(0)\|_2,   \  \  \| \Delta_{\mathbf{H}_{M}(1)} \|_2\leq \epsilon \| \mathbf{H}_{M}(1)\|_2,
\end{equation}
for some $\varepsilon>0$.
Then from (\ref{pert}) taking into account (\ref{normmat1}),  we obtain
\begin{align*}
|z^{(1)}|&=\left| \frac{\boldsymbol{p}^{T} \,  (\Delta_{\mathbf{H}_{M}(1)}-z \, \Delta_{\mathbf{H}_{M}(0)}) \,  \boldsymbol{q}}{\boldsymbol{p}^{T} \,  \mathbf{H}_{M}(0) \,  \boldsymbol{q}}  \right| \notag \\
& \leq \frac{\|\boldsymbol{p}^{T}\|_2  \,  \|\Delta_{\mathbf{H}_{M}(1) }-z \, \Delta_{\mathbf{H}_{M}(0)}\|_2 \,  \|\boldsymbol{q}\|_2}{|\boldsymbol{p}^{T}\,  \mathbf{H}_{M}(0) \,  \boldsymbol{q}|} \leq \epsilon \rho,  \label{pertL}
\end{align*}
where $\rho$ is the  \textit{sensitivity with respect to unstructured perturbations} of the eigenvalue $z$ of the GEP (\ref{gep}) defined by
\begin{equation}\label{sen}
\rho:=  \frac{\| \boldsymbol{p}^{T}  \|_2 \,  (|z| \, \| \mathbf{H}_{M}(0)\|_2 +\| \mathbf{H}_{M}(1)\|_2) \,  \| \boldsymbol{q}\|_2}{|\boldsymbol{p}^{T}  \,  \mathbf{H}_{M}(0) \,   \boldsymbol{q}|}.
\end{equation}

Next, we derive formulas for the exact computation of sensitivities  $\rho_j$ of the eigenvalues $z_j$,  $j=1,\ldots,M$ of the GEP (\ref{gep}) as well as their upper bounds.  Since the Hankel matrices $\mathbf{H}_{M}(0)$ and $\mathbf{H}_{M}(1)$ satisfy factorizations (\ref{fachan}) and (\ref{fachan1}) for $N=L=M$,
left $\boldsymbol{p}^{(j)}$ and right $\boldsymbol{q}^{(j)}$ eigenvectors  corresponding to the eigenvalue $z_j$   can be determined by formulas 
\begin{align}
\boldsymbol{q}^{(j)}=\boldsymbol{p}^{(j)}=\left( \mathbf{V}^{T}_{M}(\boldsymbol{z}) \right)^{-1} \, \boldsymbol{e}_j, \label{einvec11}
\end{align}
where $\boldsymbol{e}_j$ is the $j$th canonical column vector in $\mathbb{R}^{M}$.
Then using again the factorization (\ref{fachan})  of the  Hankel matrix $\mathbf{H}_{M}(0)$,   we get
\begin{align}
(\boldsymbol{p}^{(j)} )^{T}  \,  \mathbf{H}_{M}(0) \,   \boldsymbol{q}^{(j)} &=\left( \left( \mathbf{V}^{T}_{M} (\boldsymbol{z})  \right)^{-1} \, \boldsymbol{e}_j \right)^T \,  \mathbf{H}_{M}(0) \,  \left( \mathbf{V}^{T}_{M} (\boldsymbol{z}) \right)^{-1} \, \boldsymbol{e}_j \notag \\
& = \boldsymbol{e}_j^{T}  \,  \mathbf{\Gamma}_M \,   \boldsymbol{e}_j= \gamma_j.   \label{gam}
\end{align}
From (\ref{sen}) taking into account (\ref{gam}),  we obtain that the sensitivity $\rho_j$  is given by 
\begin{equation}\label{senhan}
\rho_j=\frac{1}{|\gamma_j|} \| (\boldsymbol{p}^{(j)})^{T} \|_2 \,  (|z_j| \, \| \mathbf{H}_{M}(0)\|_2+   \| \mathbf{H}_{M}(1)\|_2 ) \,  \| \boldsymbol{q}^{(j)}  \|_2.
\end{equation}
The ill-disposedness of the eigenvalue $z_j$ with respect to unstructured perturbations can be measured by checking whether $\rho_j$  is ''large``.
Applying  (\ref{fachan}),   (\ref{fachan1}) and (\ref{einvec11}) to (\ref{senhan}),  we obtain  the  following final formula for the computation  of   sensitivities $\rho_j$
\begin{align}
\rho_j &=\frac{1}{|\gamma_j|} \left \|  \left( \mathbf{V}^{T}_{M}(\boldsymbol{z}) \right)^{-1} \, \boldsymbol{e}_j \right \|_2^{2}   \left(|z_j| \left \| \mathbf{V}_{M} (\boldsymbol{z}) \,  \mathbf{\Gamma}_M \,  \mathbf{V}^{T}_{M} (\boldsymbol{z})  \right \|_2+   \left \| \mathbf{V}_{M} (\boldsymbol{z}) \,  \mathbf{Z}_M \,  \mathbf{\Gamma}_M \,  \mathbf{V}^{T}_{M} (\boldsymbol{z})  \right \|_2 \right).  \label{senk1}
\end{align}

We proceed with the estimation of (\ref{senhan}).  Since   (\ref{fachan}) and  (\ref{fachan1})  implies
\begin{align}
\|\mathbf{H}_{M}(0)\|_2 & \leq \|\mathbf{V}_{M} (\boldsymbol{z})\|_2  \,  \|\mathbf{\Gamma}_M\|_2 \,  \| \mathbf{V}^{T}_{M} (\boldsymbol{z}) \|_2 , \label{normmatnew} \\
\|\mathbf{H}_{M}(1)\|_2 & \leq \|\mathbf{V}_{M} (\boldsymbol{z})\|_2  \,  \|\mathbf{Z}_M \,  \mathbf{\Gamma}_M\|_2 \,  \| \mathbf{V}^{T}_{M}(\boldsymbol{z}) \|_2,   \label{normmat}
\end{align}
 and taking into account 
 the inequality  $\| \mathbf{X} \, \boldsymbol{e}_j \|_2\leq \|\mathbf{X} \|_2$,  we conclude that
$$
\rho_j\leq \frac{|z_j| \,   \| \mathbf{\Gamma}_M\|_2+\|\mathbf{Z}_M \,   \mathbf{\Gamma}_M\|_2}{|\gamma_j|} \,   \kappa^{2} \left( \mathbf{V}_{M} (\boldsymbol{z}) \right) ,
$$
where $\kappa( \mathbf{X}):=\|\mathbf{X}\|_2 \, \|\mathbf{X}^{-1}\|_2$ denotes the condition number of the matrix $\mathbf{X}$.  For diagonal matrices $\mathbf{\Gamma}_M$ and $\mathbf{Z}_M \,   \mathbf{\Gamma}_M$,  it holds
$$
 \|\mathbf{\Gamma}_M\|_2=\max\limits_{k=1,\ldots,M}|\gamma_k|, \ \ \  \|\mathbf{Z}_M \,   \mathbf{\Gamma}_M\|_2=\max\limits_{k=1,\ldots,M} |z_k \gamma_k|.
$$
Thus, we derive the following upper bounds for sensitivities $\rho_j$,
\begin{equation}\label{senhan1}
\rho_j \leq \zeta_j \,   \kappa^{2} \left( \mathbf{V}_{M} (\boldsymbol{z}) \right) ,
\end{equation}
where $\zeta_j$ is defined by 
\begin{equation}\label{zeta11}
\zeta_j:
=\frac{1}{|\gamma_j|} \left(|z_j| \max\limits_{k=1,\ldots,M}  |\gamma_k| + \max\limits_{k=1,\ldots,M}  |\gamma_k z_k|  \right) .
\end{equation}
 Note that we can  estimate (\ref{senhan1}) further by applying known estimates for condition numbers of  Vandermonde matrices \cite{P16, AB2019}.

Defining $\boldsymbol{\zeta}:=(\zeta_1,\ldots,\zeta_M)^T \in \rr^{M}$  and the vector of  sensitivities  $\boldsymbol{\rho}:=(\rho_1,\ldots,\rho_M)^{T}\in \rr^{M}$,  and  applying (\ref{senhan1}), we get for the $\ell_2$-norm of the vector  $\boldsymbol{\rho}$,
\begin{equation*}\label{l2sen}
\|\boldsymbol{\rho}\|_2 \leq \|\boldsymbol{\zeta}\|_2 \,   \kappa^{2} \left( \mathbf{V}_{M} (\boldsymbol{z}) \right)  .
\end{equation*}
From (\ref{senk1}),  taking into account  (\ref{normmatnew}), (\ref{normmat}) and (\ref{zeta11}),   we obtain 
$$
\rho_j\leq \zeta_j  \bigg \| \left( \mathbf{V}^{T}_{M}(\boldsymbol{z}) \right)^{-1} \, \boldsymbol{e}_j  \bigg \|_2^{2} \,  \big \|\mathbf{V}_{M} (\boldsymbol{z}) \big \|_2^{2} .
$$
By summing up the last inequalities  over $j=1,\ldots,M$,  we arrive at
$$
\sum\limits_{j=1}^{M} \rho_j \leq \zeta_{\max}   \big \|\mathbf{V}_{M} (\boldsymbol{z}) \big \|_2^{2} \,  \sum\limits_{j=1}^{M}  \bigg\| \left( \mathbf{V}^{T}_{M}(\boldsymbol{z}) \right)^{-1} \, \boldsymbol{e}_j  \bigg\|_2^{2}, 
$$
where $\zeta_{\max}:=\max_{j=1,\ldots,M} \zeta_j$.
Finally, using the property $\sum_{j=1}^{M} \|\mathbf{X} \boldsymbol{e}_j \|_2^{2}=  \|\mathbf{X} \|_{\mathrm{F}}^{2}$,  we derive the following upper bound
\begin{equation*}\label{l1n}
\| \boldsymbol{\rho}\|_1 \leq \zeta_{\max}   \big \|\mathbf{V}_{M} (\boldsymbol{z}) \big \|_2^{2} \,   \big \| \mathbf{V}^{-1} _{M}(\boldsymbol{z})  \big  \|_{\mathrm{F}}^{2} 
\end{equation*}
 for the $\ell_1$-norm of the sensitivity vector $\boldsymbol{\rho}$.

\subsection{Sensitivities of Hankel  pencil eigenvalues with structured perturbations}
\label{stsen}

In this subsection, we consider structured perturbations of the GEP (\ref{gep}) for the matrix pencil (\ref{mp11}). Now  the perturbation matrices $\Delta_{\mathbf{H}_{M}(0)}$ and $\Delta_{\mathbf{H}_{M}(1)}$ in (\ref{pertsys}) have a particular structure,  they are again Hankel matrices.  We give a short description of the results from \cite{Zhphd}.  Let $\boldsymbol{\delta}_{2M}:=(\delta_0,\delta_1,\ldots,\delta_{2M-1})^{T}$  be a vector of perturbations and let
$$
\Delta_{\mathbf{H}_{M}(0)}= (\delta_{k+\ell})_{k,\ell=0}^{M-1,M-1} ,  \  \Delta_{\mathbf{H}_{M}(1)}= (\delta_{k+\ell+1})_{k,\ell=0}^{M-1,M-1}.
$$
For the first order approximation of the eigenvalue perturbation $z_j^{(1)}$,  $j=1,\ldots,M$ from (\ref{pert}) using (\ref{gam}),  we obtain 
\begin{equation}\label{pertt}
z^{(1)}_j=\frac{1}{\gamma_j} \, ( \boldsymbol{p}^{(j)} )^{T} \,  (\Delta_{\mathbf{H}_{M}(1)}-z_j \, \Delta_{\mathbf{H}_{M}(0)}) \,  \boldsymbol{q}^{(j)}.
\end{equation}
Define further  matrices  $\mathbf{Q}^{(j)}_{M,2M}(0):=\left(Q^{(j)}_{k,m}(0)\right)_{k,m=1}^{M,2M}$ and $\mathbf{Q}^{(j)}_{M,2M}(1):=\left(Q^{(j)}_{k,m}(1)\right)_{k,m=1}^{M,2M}$ by their entries
$$
Q^{(j)}_{k,m}(0):=
\begin{cases}
q^{(j)}_{m-k+1},   & 1 \leq m-k+1\leq M, \\
0,   &  \text{otherwise},
\end{cases}
$$
and 
$$
Q^{(j)}_{k,m}(1):=
\begin{cases}
q^{(j)}_{m-k},   & 1 \leq m-k\leq M, \\
0,   &  \text{otherwise}. 
\end{cases}
$$
Then
$$
\Delta_{\mathbf{H}_{M}(0)} \,  \boldsymbol{q}^{(j)}= \mathbf{Q}_{M,2M}^{(j)}(0) \, \boldsymbol{\delta}_{2M},  \  \Delta_{\mathbf{H}_{M}(1)} \, \boldsymbol{q}^{(j)}= \mathbf{Q}_{M,2M}^{(j)}(1) \, \boldsymbol{\delta}_{2M}.
$$
Since for $m=1,\ldots,2M$
\begin{align*}
\left( (\boldsymbol{p}^{(j)})^T \,    \mathbf{Q}^{(j)}_{M,2M}(0)   \right)_m & =\sum\limits_{k=1}^{2M} p^{(j)}_{k} \,  Q^{(j)}_{k,m}(0) = \sum\limits_{k=\max\{1,m-M+1\}}^{\min\{M,m \}} p^{(j)}_{k} \,  q^{(j)}_{m-k+1},   \\
\left( (\boldsymbol{p}^{(j)})^T \,    \mathbf{Q}_{M,2M}^{(j)}(1)   \right)_m & =\sum\limits_{k=1}^{M} p^{(j)}_{k} \,  Q_{k,m}^{(j)}(1)= \sum\limits_{k=\max\{1,m-M\}}^{\min\{M,m-1 \}} p^{(j)}_{k} \,  q^{(j)}_{m-k},
\end{align*}
and taking into account (\ref{pertt}),  we
define  vectors $\boldsymbol{\mathcal{S}}^{(j)}_{2M}:=\left(S^{(j)}_1, \ldots, S^{(j)}_{2M}\right)^{T} \in \cc^{2M}$,  $j=1,\ldots,M$ by their components
\begin{equation*}\label{vecs}
S^{(j)}_m:=\frac{1}{\gamma_j}\left(  \sum\limits_{k=\max\{1,m-M\}}^{\min\{M,m-1 \}} p^{(j)}_{k} \,  q^{(j)}_{m-k} - z_j \,  \sum\limits_{k=\max\{1,m-M+1\}}^{\min\{M,m \}} p^{(j)}_{k} \,  q^{(j)}_{m-k+1} \right),
\end{equation*}
for $m=1,\ldots,2M$.  Then  $z_j^{(1)}$ as in (\ref{pertt}) can be presented by
\begin{equation}\label{pe}
z_j^{(1)}= (\boldsymbol{\mathcal{S}}^{(j)}_{2M})^T \,  \boldsymbol{\delta}_{2M}.
\end{equation}
The component $S^{(j)}_m$ shows the influence of the perturbations of the $m$th measurement (in this case it is the function value $f(m-1)$)  on the $j$th eigenvalue $z_j$ of the Hankel matrix pencil (\ref{mp11}).  Using (\ref{einvec11}),  we rewrite  components $S^{(j)}_m$  as follows
\begin{align*}
S^{(j)}_m =& \frac{1}{\gamma_j} \bigg( \,  \sum\limits_{k=\max\{1,m-M\}}^{\min\{M,m-1 \}} \bigg({\left( \mathbf{V}^{T}_{M} (\boldsymbol{z}) \right)^{-1}}\bigg)_{k,j} \,  \bigg( \left( \mathbf{V}^{T}_{M} (\boldsymbol{z})\right)^{-1} \bigg)_{m-k,j} \notag \\
& - z_j \,  \sum\limits_{k=\max\{1,m-M+1\}}^{\min\{M,m \}} \bigg( \left( \mathbf{V}^{T}_{M} (\boldsymbol{z}) \right)^{-1} \bigg)_{k,j} \, \bigg(  \left( \mathbf{V}^{T}_{M} (\boldsymbol{z}) \right)^{-1} \bigg)_{m-k+1,j} \bigg). \label{sent}
\end{align*}
Formula (\ref{pe}) can be also written in a matrix form 
\begin{equation*}\label{pe1}
\boldsymbol{z}^{(1)}= \mathbf{\mathcal{S}}_{M,2M} \, \boldsymbol{\delta}_{2M},
\end{equation*} 
where $\boldsymbol{z}^{(1)}:=(z_1^{(1)},\ldots, z_M^{(1)})^{T}$ and  $\mathbf{\mathcal{S}}_{M,2M}:= (S_{j,m})_{j,m=1}^{M,2M} \in \cc^{M\times 2M}$ be a matrix  with the entries $S_{j,m}:=S^{(j)}_m$, which quantify the perturbation of the eigenvalues of the matrix pencil (\ref{mp11}).

Finally,  we define  the  sensitivity $\eta_{j,m}$ of the eigenvalue $z_j$  of the matrix pencil (\ref{mp11}) with respect to the $m$th measurement by
\begin{equation*}
\eta_{j,m}:=|S_{j,m}|,
\end{equation*}
for $j=1,\ldots,M$ and $m=1,\ldots,2M$. By $\boldsymbol{\eta}_{2M}^{(j)}:=(\eta_{j,1},\ldots, \eta_{j,2M})^T \in \mathbb{R}^{2M}$ we denote the corresponding vector of sensitivities for the eigenvalue $z_j$.  Then
\textit{the sensitivity $\eta_j$  with respect to structured perturbations}  for the eigenvalue $z_j$  of the Hankel matrix pencil (\ref{mp11}) is defined as the  $\ell_2$-norm of the vector  $\boldsymbol{\mathcal{S}}^{(j)}_{2M} \in \cc^{2M}$.  Thus, we further define  a vector $\boldsymbol{\eta}:=(\eta_1,\ldots, \eta_M)^{T} \in \mathbb{R}^{M}$ by its components
\begin{equation*}\label{stp}
\eta_j:=\|\boldsymbol{\mathcal{S}}^{(j)}_{2M}\|_2=\sqrt{\eta_{j,1}^{2}+\ldots+\eta^{2}_{j,2M}},   \quad  j=1,\ldots,M.
\end{equation*}
The ill-disposedness of the eigenvalue $z_j$ with respect to structured perturbations can be measured by checking whether $\eta_j$  is ''large``.
The entry $\eta_j$ determines also the standard deviation of  the eigenvalue perturbation $z_j^{(1)}$ (see  \cite{Zhphd} for details).

\section{Recovery of rational  functions via Hankel matrix pencil method and sensitivities of the poles}
\label{method}

In this section,  we describe the main ideas of our new method for the recovery of  rational functions (\ref{rat}) and present the sensitivity analysis for the poles.

\subsection{Description of the method}
\label{method1}

Let $r(t)$ be a rational function defined by (\ref{rat}) and
by $\mathbb{D}:=\{z \in \cc:  \, |z|=1 \}$ we denote the unit circle in the complex plane.  
Following approach in \cite{Y22},  for $k \in \mathbb{Z}\setminus\{0\}$  we consider the integral 
$$
\frac{1}{2\pi \mathrm{i}} \int\limits_{\mathbb{D}} \frac{r(t)}{t^{k}} \, \frac{\mathrm{d}t}{t},
$$
which is closely related to the classical Fourier coefficient $\widehat{r}_k$ of the function $r(\mathrm{e}^{\mathrm{i}\varphi})$, 
\begin{equation} \label{rkf}
\widehat{r}_k:=\frac{1}{2\pi} \int\limits_0^{2\pi} r(\mathrm{e}^{\mathrm{i}\varphi} ) \mathrm{e}^{-\mathrm{i}k \varphi} \, \mathrm{d} \varphi = \frac{1}{2\pi \mathrm{i}} \int\limits_{\mathbb{D}} \frac{r(t)}{t^{k}} \, \frac{\mathrm{d}t}{t}. 
\end{equation}
Then (see,  for example,  \cite{Y22})
\begin{equation}\label{coef}
\widehat{r}_k=
\begin{cases}
\sum\limits_{|z_j|<1} \gamma_j z_j^{-(k+1)},  & k = -1, -2,\ldots , \\
-\sum\limits_{|z_j|>1} \gamma_j z_j^{-(k+1)},  & k= 1, 2 , \ldots. 
\end{cases}
\end{equation}
From (\ref{coef}),  we conclude that Fourier coefficients $\widehat{r}_k$ with negative indices $k=-1,-2,  \ldots$ contain information about poles located inside  the unite circle,  and Fourier coefficients $\widehat{r}_k$ with positive indices $k=1,2,\ldots$ contain information about poles located outside  the unite circle.

Taking into account (\ref{coef}),   we define an exponential sum $f(t)$ by 
\begin{equation}\label{exp}
f(t)=
\begin{cases}
\sum\limits_{|z_j|<1} \gamma_j z_j^{t},  & t \geq 0, \\
-\sum\limits_{|z_j|>1} \gamma_j z_j^{t},  & t < 0, 
\end{cases}
\end{equation}
and
from (\ref{coef}) we obtain the following interpolation property
\begin{equation}\label{int}
f(-(k+1))=\widehat{r}_k,   \quad  k \in \zz\setminus\{0\}.
\end{equation}
Thus, the rational reconstruction problem (\ref{rat}) can be reformulated as the exponential interpolation problem (\ref{int}).  Our further goal is to reconstruct the exponential sum (\ref{exp}),  using samples $f(m)$,  $m\in \zz \setminus \{-1\}$ as inputs.  From (\ref{exp}) and (\ref{int}),  we conclude that values $f(m)$ with $m=0,1,\ldots$ contain information about nodes $z_j$ located inside $\mathbb{D}$,  and values $f(m)$ with $m=-2,-3,\ldots$ contain information about nodes $z_j$  located outside $\mathbb{D}$. 

Let $M_1$ and $M_2$ be the numbers of nodes $z_j$ of $f$ as in (\ref{exp}) located inside $\mathbb{D}$ and outside  $\mathbb{D}$,  respectively.    
Thus, we can rewrite the exponential sum (\ref{exp}) as follows
\begin{equation}\label{exp1}
f(t)=
\begin{cases}
\sum\limits_{n=1}^{M_1} \gamma^{(-)}_{n} (z^{(-)}_{n})^{t},  & t \geq 0, \\
-\sum\limits_{n=1}^{M_2} \gamma^{(+)}_{n} (z^{(+)}_{n})^{t}  & t <0, 
\end{cases}
\end{equation}
i.e.  by $z^{(-)}_{n}$,  $n=1,\ldots,M_1$, and $z^{(+)}_{n}$,  $n=1,\ldots,M_2$, we denote nodes located inside $\mathbb{D}$ and outside $\mathbb{D}$, respectively, and let $\gamma^{(-)}_{n} $,  $n=1,\ldots,M_1$, and $\gamma^{(+)}_{n} $,  $n=1,\ldots,M_2$, be the corresponding coefficients.

Next, we present some modification of the ESPRIT method to solve  our special parameter estimation problem (\ref{exp1}).  Note that  ESPRIT  is not the only method to apply in the frequency domain.  We can use also the other methods for the exponential recovery,  for example,  MPM (matrix-pencil method) \cite{PT2013}, approximate Prony method \cite{PT10}, or AAK theory \cite{PP} (as in   \cite{Y22}).

By $\boldsymbol{z}^{(-)} \in \mathbb{C}^{M_1}$ and $\boldsymbol{z}^{(+)} \in \mathbb{C}^{M_2}$ we denote two vectors
$$
\boldsymbol{z}^{(-)}:=(z^{(-)}_{1},\ldots, z^{(-)}_{M_1})^{T} ,  \  \ \boldsymbol{z}^{(+)}:=\left(\frac{1}{z_{1}^{(+)}},\ldots,  \frac{1}{z^{(+)}_{M_2}}\right)^{T},
$$
and let us also define the following  diagonal matrices 
$$
\mathbf{\Gamma}^{(-)}_{M_1}:=\mathrm{diag} (\gamma^{(-)}_{n})_{n=1}^{M_1},   \  \mathbf{\Gamma}^{(+)}_{M_2}:=\mathrm{diag} (\gamma^{(+)}_{n})_{n=1}^{M_2},  \  \mathbf{Z}_{M_1}^{(+)}:=\mathrm{diag}\left( z^{(-)}_{n} \right)_{n=1}^{M_1},  \ \mathbf{Z}_{M_2}^{(+)}:=\mathrm{diag}\left( \frac{1}{z^{(+)}_{n}} \right)_{n=1}^{M_2}.
$$ 
Let further $L,  N \in \nn$ be given such that $M_1,M_2 \leq L \leq N$, and  $L$ is the upper bound for $M$.   We use values $f(k)$,  $k=0,\ldots,2N-1$ ($\widehat{r}_{-k}$,  $k=1, \ldots, 2N$) for recovering nodes $z^{(-)}_{n}$,  $n=1,\ldots,M_1$ located inside $\mathbb{D}$,  and values  $f(-k)$,  $k=2,\ldots,2N+1$ ($\widehat{r}_{k}$,  $k=1, \ldots, 2N$) for recovering nodes $z^{(+)}_{n}$,  $n=1,\ldots,M_2$ located outside $\mathbb{D}$.
According to the ESPRIT method described in Section \ref{esprit},  to reconstruct nodes $z^{(-)}_{n}$,  $n=1,\ldots,M_1$, and  $z^{(+)}_{n}$,  $n=1,\ldots,M_2$,  we employ the 
 Hankel matrices
 \begin{align*}
 \mathbf{H}^{(-)}_{2N-L,L+1}:=\left(f(\ell+k)\right)_{k,\ell=0}^{2N-L-1,L}= \left(\widehat{r}_{-(k+\ell)}\right)_{k=0,\ell=1}^{2N-L-1,L+1}, 
 \end{align*}
and
 $$
  \mathbf{H}^{(+)}_{2N-L,L+1}:=(f(-(\ell+k+2)))_{k,\ell=0}^{2N-L-1,L}= \left(\widehat{r}_{k+\ell}\right)_{k=0,\ell=1}^{2N-L-1,L+1},
 $$
 respectively.
We consider factorizations of these matrices
\begin{align}
 \mathbf{H}^{(-)}_{2N-L,L+1}&=\mathbf{V}_{2N-L,M_1} (\boldsymbol{z}^{(-)}) \,  \mathbf{\Gamma}^{(-)}_{M_1} \,  \mathbf{V}^{T}_{L+1,M_1} (\boldsymbol{z}^{(-)}), \label{fac1} \\
  \mathbf{H}^{(+)}_{2N-L,L+1}&=-\mathbf{V}_{2N-L,M_2} (\boldsymbol{z}^{(+)}) \,  \mathbf{\Gamma}^{(+)}_{M_2} \,  \left(\mathbf{Z}_{M_2}^{(+)}\right)^{2}  \,  \mathbf{V}^{T}_{L+1,M_2} (\boldsymbol{z}^{(+)}). \label{fac2}
\end{align} 
  Since all Vandermonde  matrices in (\ref{fac1}) and (\ref{fac2})  have full column ranks,   $\mathbf{H}_{2N-L,L+1}^{(-)}$ has rank $M_1$ and $\mathbf{H}_{2N-L,L+1}^{(+)}$ has rank $M_2$.  Therefore, in case of exact given data,  the numbers $M_1$ and $M_2$ are reconstructed as numerical ranks of matrices $   \mathbf{H}_{2N-L,L+1}^{(-)}$ and $\mathbf{H}_{2N-L,L+1}^{(+)}$,  respectively.

Let us start with the case of nodes  $z^{(-)}_{n}$,  $n=1,\ldots,M_1$,   located inside   $\mathbb{D}$.  Since  the Hankel matrix $ \mathbf{H}^{(-)}_{2N-L,L+1}$ have  the structure similar to those  of the Hankel matrix $\mathbf{H}_{2N-L,L+1}$ as in (\ref{factH}),  following ideas of  Section \ref{esprit},  we conclude that 
nodes $z^{(-)}_{n}$,  $n=1,\ldots,M_1$ can be computed as eigenvalues of the Hankel matrix pencil
\begin{equation}\label{matp1}
z \mathbf{H}^{(-)}_{2N-L,L}(0) -  \mathbf{H}^{(-)}_{2N-L,L}(1),
\end{equation}
where for $s=0,1$,
$$
\mathbf{H}^{(-)}_{2N-L,L}(s):=\mathbf{H}^{(-)}_{2N-L,L+1}(1:2N-L,1+s:L+s)=(f(\ell+k+s))_{k,\ell=0}^{2N-L-1,L-1}.
$$

Consider now the case of nodes  $z^{(+)}_{n}$,  $n=1,\ldots,M_2$,   located outside  $\mathbb{D}$.  For $s=0,1$,  define matrices 
\begin{align*}
\mathbf{H}^{(+)}_{2N-L,L}(s):=\mathbf{H}^{(+)}_{2N-L,L+1}(1:2N-L,1+s:L+s)=(f(-(\ell+k+s+2)))_{k,\ell=0}^{2N-L-1,L-1},
\end{align*}
and write their factorizations
\begin{align}
 \mathbf{H}^{(+)}_{2N-L,L}(0)&=-\mathbf{V}_{2N-L,M_2} (\boldsymbol{z}^{(+)}) \,  \mathbf{\Gamma}^{(+)}_{M_2} \,  \left(\mathbf{Z}_{M_2}^{(+)}\right)^{2}  \,  \mathbf{V}^{T}_{L,M_2} (\boldsymbol{z}^{(+)}),   \label{fac21n}  \\
  \mathbf{H}^{(+)}_{2N-L,L}(1) &=-\mathbf{V}_{2N-L,M_2} (\boldsymbol{z}^{(+)}) \,  \mathbf{\Gamma}^{(+)}_{M_2} \,  \left(\mathbf{Z}_{M_2}^{(+)}\right)^{3}  \,  \mathbf{V}^{T}_{L,M_2} (\boldsymbol{z}^{(+)}). \label{fac22n}
\end{align}
According to (\ref{fac21n}) and (\ref{fac22n}) we have
 $$ z \,  {\mathbf H}^{(+)}_{2N-L,L}(0) - {\mathbf H}^{(+)}_{2N-L,L}(1) = -\mathbf{V}_{2N-L,M_2} (\boldsymbol{z}^{(+)})   \, \mathrm{diag}  \left( z - \frac{1}{z^{(+)}_{n}} \right)_{n=1}^{M_2}  \,   \mathbf{\Gamma}^{(+)}_{M_2} \,  \left(\mathbf{Z}_{M_2}^{(+)}\right)^{2} \,  \mathbf{V}^{T}_{L,M_2} (\boldsymbol{z}^{(+)}),$$
Thus,  the inverse $\frac{1}{z^{(+)}_{n}}$ of nodes $z^{(+)}_{n}$ for $n=1,\ldots,M_2$ are computed as eigenvalues of the Hankel matrix pencil 
\begin{equation}\label{matp2}
z \mathbf{H}^{(+)}_{2N-L,L}(0) -  \mathbf{H}^{(+)}_{2N-L,L}(1).
\end{equation}

If we assume that $M_1$ and $M_2$ are known,  we can reconstruct
nodes $z^{(-)}_{n}$,  $n=1,\ldots,M_1$ as eigenvalues of the Hankel matrix pencil ($N=L=M_1$ in (\ref{matp1}))
\begin{equation}\label{mp1}
z \mathbf{H}^{(-)}_{M_1}(0) -  \mathbf{H}^{(-)}_{M_1}(1),
\end{equation}
and we can reconstruct  the inverse $\frac{1}{z^{(+)}_{n}}$ of nodes $z^{(+)}_{n}$, $n=1,\ldots,M_2$, as eigenvalues of the Hankel matrix pencil ($N=L=M_2$ in (\ref{matp2}))
\begin{equation}\label{mp2}
z \mathbf{H}^{(+)}_{M_2}(0) -  \mathbf{H}^{(+)}_{M_2}(1).
\end{equation}

Having nodes $z^{(-)}_{n}$  and $z^{(+)}_{n}$,  coefficients $\gamma^{(-)}_{n}$ and $\gamma^{(+)}_{n}$ can be computed as  least squares solutions of the overdetermined systems of linear equations
$$
\sum\limits_{n=1}^{M_1} \gamma^{(-)}_{n} (z^{(-)}_{n})^{-(k+1)} =\widehat{r}_{-k},  \quad k=1,\ldots,2N,
$$
and
$$
-\sum\limits_{n=1}^{M_2} \gamma^{(+)}_{n} (z^{(+)}_{n})^{-(k+1)}= \widehat{r}_k,  \quad k=1,\ldots,2N,
$$
respectively.

 \begin{remark}
 \label{rem41}
\textbf{ 1.} Note that Fourier coefficients $\hat{r}_k$ can also be computed applying quadrature formulas on some grid of points.  In that case,  we  use the function values of $r$ instead of its Fourier coefficients as inputs for the algorithm.   Let us chose points $t_n=\mathrm{e}^{\frac{2\pi \mathrm{i} n}{4N}} $, $n=0,...,4N-1$ and let $(R_k)_{k=0}^{4N-1}$ be the DFT of $(r(t_n))_{n=0}^{4N-1}$. Then $\hat{r}_k$ in (\ref{rkf}) can be determined by $\hat{r}_k=R_{k+2N}$ for $k=-2N,...,2N-1$.
 This question was discussed  in details in \cite{Y22},  where the author considered not only the points on the unit circle, but also the Matsubara model as a set of grid points for the evaluation of the transform $\widehat{r}_k$.  On the one hand,  we are interested in the presentation of a modification of the method from \cite{Y22} and in the explanation of the results of this method for the recovery from noisy data with the use of the sensitivities. In that case, we have to employ the given noisy coefficients $\widehat{r}_k$ as input data. On the other hand, we also want to compare the performance of known methods for rational approximation with our new approach based on Hankel pencil. To guarantee fair comparison, we assume to have given noisy values $(r(t_n))_{n=0}^{4N-1}$.

\textbf{ 2.} The method based on using the given coefficients $\hat{r}_k$, $k=-2N,...,2N$ has the overall computational cost of $\mathcal{O}(N L^{2})$ flops, which is determined by the complexity of the computation of SVDs of Hankel matrices $ \mathbf{H}^{(-)}_{2N-L,L+1}$ and $ \mathbf{H}^{(+)}_{2N-L,L+1}$. 
The complexity of the modified algorithm based on employing the values $(r(t_n))_{n=0}^{4N-1}$ is $\mathcal{O}(N (L^{2}+\log N))$, since the additional step of the computation of FFT of $(r(t_n))_{n=0}^{4N-1}$ requires  $\mathcal{O}(N \log N)$ operations. 
 \end{remark}

\subsection{Sensitivities of poles }
\label{sub42}

In this subsection,  we present  sensitivity analysis for poles $z_j$ of the rational function (\ref{rat}), reconstructed with the method described in Subsection \ref{method1}.
First, we consider the case of unstructured perturbations.  Poles $z_{n}^{(-)}$,  $n=1,\ldots,M_1$ located inside the unit circle $\mathbb{D}$  are reconstructed as eigenvalues of the matrix pencil  (\ref{mp1}).  Since Hankel matrices in (\ref{mp1}) have  structures similar to those  of Hankel matrices in  (\ref{mp11}),  we  apply the ideas from Subsection \ref{unsen}  without any changes and obtain  the formula for the exact computation of  sensitivities $\rho_n^{(-)}$ with respect to unstructured perturbations of the poles $z^{(-)}_{n}$,  $n=1,\ldots,M_1$, 
\begin{align}
\rho_n^{(-)} &=\frac{1}{\left|\gamma_{n}^{(-)}\right|} \left \| \left( \mathbf{V}^{T}_{M_1}(\boldsymbol{z}^{(-)}) \right)^{-1} \, \boldsymbol{e}_n  \right \|_2^{2} \notag \\
& \times \left(|z_{n}^{(-)}| \left \| \mathbf{V}_{M_1} (\boldsymbol{z}^{(-)}) \,  \mathbf{\Gamma}_{M_1} \,  \mathbf{V}^{T}_{M_1} (\boldsymbol{z}^{(-)})  \right \|_2+   \left \| \mathbf{V}_{M_1} (\boldsymbol{z}^{(-)}) \,  \mathbf{Z}_{M_1} \,  \mathbf{\Gamma}_{M_1} \,  \mathbf{V}^{T}_{M_1} (\boldsymbol{z}^{(-)})  \right \|_2 \right)  , \label{senk}
\end{align}
where $\boldsymbol{e}_n$  is the $n$th canonical column vector in $\mathbb{R}^{M_1}$. Correspondingly, 
for the upper bound of $\rho_n^{(-)}$ we have
\begin{equation}\label{se1}
\rho_n^{(-)} \leq \zeta_n^{(-)} \,   \kappa^{2} \left( \mathbf{V}_{M_1} (\boldsymbol{z}^{(-)})\right) ,
\end{equation}  
with
$$
\zeta_n^{(-)}:=\frac{1}{\left|\gamma_{n}^{(-)}\right|} \left(\left|z_{n}^{(-)}\right| \max\limits_{k=1,\ldots,M_1}  \left|\gamma_{k}^{(-)}\right| + \max\limits_{k=1,\ldots,M_1}  \left|\gamma_{k}^{(-)} z_{k}^{(-)}\right|  \right) .
$$

Now consider the case of  sensitivities $\rho_n^{(+)}$ with respect to unstructured perturbations for the poles $z^{(+)}_{n}$,  $n=1,\ldots,M_2$ located  outside $\mathbb{D}$. As was already mentioned,  the inverse  $\frac{1}{z^{(+)}_{n}}$ of the poles $z^{(+)}_{n}$ are reconstructed as eigenvalues of
 the matrix pencil (\ref{mp2}).   Thus, we compute sensitivities $ \rho_n^{(+)}$  of the inverse $\frac{1}{z^{(+)}_{n}}$  of the poles $z^{(+)}_{n}$,  but we  say that $ \rho_n^{(+)}$ are also sensitivities of $z^{(+)}_{n}$. 
 
 To present formulas for sensitivities  $\rho_n^{(+)}$,  
  we have to apply some modifications to the results from  Subsection \ref{unsen}.  Matrices  $ \mathbf{H}^{(+)}_{M_2}(0)$ and $ \mathbf{H}^{(+)}_{M_2}(1)$ satisfy factorizations
  \begin{equation}\label{fac21}
 \mathbf{H}^{(+)}_{M_2}(0)=-\mathbf{V}_{M_2} (\boldsymbol{z}^{(+)}) \,  \mathbf{\Gamma}^{(+)}_{M_2} \,  \left(\mathbf{Z}_{M_2}^{(+)}\right)^{2}  \,  \mathbf{V}^{T}_{M_2} (\boldsymbol{z}^{(+)}),
 \end{equation}
   \begin{equation}\label{fac22}
 \mathbf{H}^{(+)}_{M_2}(1)=-\mathbf{V}_{M_2} (\boldsymbol{z}^{(+)}) \,  \mathbf{\Gamma}^{(+)}_{M_2} \,  \left(\mathbf{Z}_{M_2}^{(+)}\right)^{3}  \,  \mathbf{V}^{T}_{M_2} (\boldsymbol{z}^{(+)}).
 \end{equation}
Thus,   right and left eigenvectors of the   GEP of the Hankel matrix pencil (\ref{mp2}) corresponding to the eigenvalue $\frac{1}{z^{(+)}_{n}}$  are given by formulas 
\begin{align}
\boldsymbol{q}^{(n)}= \left( \mathbf{V}^{T}_{M_2}(\boldsymbol{z}^{(+)}) \right)^{-1} \, \boldsymbol{e}_n,  \quad 
\boldsymbol{p}^{(n)}=\left( \mathbf{V}^{T}_{M_2} (\boldsymbol{z}^{(+)})  \right)^{-1} \, \boldsymbol{e}_n, \label{einvec1}
\end{align}
respectively, where $\boldsymbol{e}_n$ is the $n$th canonical column vector in $\mathbb{R}^{M_2}$ .
Taking into account (\ref{fac21}) and (\ref{einvec1}),  we get
\begin{align}
(\boldsymbol{p}^{(n)})^{T}  \,  \mathbf{H}^{(+)}_{M_2}(0) \,   \boldsymbol{q}^{(n)} &=\left( \left( \mathbf{V}^{T}_{M_2} (\boldsymbol{z}^{(+)}) \right)^{-1} \, \boldsymbol{e}_n \right)^T \,  \mathbf{H}^{(+)}_{M_2}(0) \,   \left( \mathbf{V}^{T}_{M_2} (\boldsymbol{z}^{(+)}) \right)^{-1} \, \boldsymbol{e}_n \notag \\
& = - \boldsymbol{e}_n^{T}  \,  \mathbf{\Gamma}^{(+)}_{M_2} \,  \left(\mathbf{Z}_{M_2}^{(+)}\right)^{2} \,   \boldsymbol{e}_n=- \gamma_{n}^{(+)} \, (z_{n}^{(+)})^{-2}.   \label{gam1}
\end{align} 
 From (\ref{sen}) applying (\ref{fac21})--(\ref{gam1}),  we obtain  the formula for the exact computation  of sensitivities $\rho_n^{(+)}$ of the poles $z^{(+)}_{n}$ (the eigenvalues $\frac{1}{z^{(+)}_{n}}$ of the matrix pencil (\ref{mp2})) for $n=1,\ldots,M_2$,
\begin{align}
& \rho_n^{(+)} =\frac{1}{|\gamma_{n}^{(+)} \, (z_{n}^{(+)})^{-2}|} \bigg \|   \left( \mathbf{V}^{T}_{M_2} (\boldsymbol{z}^{(+)}) \right)^{-1} \, \boldsymbol{e}_n  \bigg \|_2^{2} \notag \\
& \times \bigg(  \frac{1}{|z_{n}^{(+)}| }   \left \| \mathbf{V}_{M_2} (\boldsymbol{z}^{(+)}) \,  \mathbf{\Gamma}^{(+)}_{M_2} \,  (\mathbf{Z}_{M_2}^{(+)})^{2}  \,  \mathbf{V}^{T}_{M_2} (\boldsymbol{z}^{(+)}) \right \|_2+   \left \| \mathbf{V}_{M_2} (\boldsymbol{z}^{(+)}) \,  \mathbf{\Gamma}^{(+)}_{M_2} \,  (\mathbf{Z}_{M_2}^{(+)})^{3}  \,  \mathbf{V}^{T}_{M_2} (\boldsymbol{z}^{(+)}) \right \|_2 \bigg) , \label{senk2}
\end{align}
and the corresponding upper bound
 $$
\rho^{(+)}_n\leq \frac{\frac{1}{|z_{n}^{(+)}| }  \,  \big \|\big(\mathbf{Z}^{(+)}_{M_2} \big)^2 \mathbf{\Gamma}_{M_2}^{(+)}\big \|_2+\big\|\big(\mathbf{Z}^{(+)}_{M_2} \big)^3 \,   \mathbf{\Gamma}^{(+)}_{M_2}\big\|_2}{|\gamma_{n}^{(+)} \, (z^{(+)}_{n})^{-2}|} \,   \kappa^{2} \left( \mathbf{V}_{M_2} (\boldsymbol{z}^{(+)}) \right) .
$$
Since
$$
 \big\|\big(\mathbf{Z}^{(+)}_{M_2} \big)^2 \mathbf{\Gamma}_{M_2}^{(+)} \big\|_2=\max\limits_{k=1,\ldots,M_2} |\gamma^{(+)}_{k} \,  (z^{(+)}_{k})^{-2}|, \ \ \ \big\|\big(\mathbf{Z}^{(+)}_{M_2} \big)^3 \,   \mathbf{\Gamma}^{(+)}_{M_2}\big\|_2=\max\limits_{k=1,\ldots,M_2} |\gamma^{(+)}_{k} \,  (z^{(+)}_{k})^{-3}|,
$$
the last inequality implies
\begin{equation}\label{se2}
\rho_n^{(+)} \leq \zeta_n^{(+)} \,   \kappa^{2} \left( \mathbf{V}_{M_2} (\boldsymbol{z}^{(+)}) \right),
\end{equation}  
where
$$
\zeta_n^{(+)}:=\frac{1}{|\gamma_{n}^{(+)} \, (z^{(+)}_{n})^{-2}|} \left(\frac{1}{|z_{n}^{(+)}| }  \max\limits_{k=1,\ldots,M_2}|\gamma^{(+)}_{k} \,  (z^{(+)}_{k})^{-2}| + \max\limits_{k=1,\ldots,M_2} |\gamma^{(+)}_{k} \,  (z^{(+)}_{k})^{-3}|  \right).
$$

Define further two vectors of  sensitivities  $\boldsymbol{\rho}^{(-)}:=(\rho_1^{(-)},\ldots,\rho^{(-)}_{M_1})^{T}\in \rr^{M_1}$ and $\boldsymbol{\rho}^{(+)}:=(\rho_1^{(+)},\ldots,\rho^{(+)}_{M_2})^{T}\in \rr^{M_2}$,  and let 
 $\boldsymbol{\zeta}^{(-)}:=(\zeta_1^{(-)},\ldots,\zeta^{(-)}_{M_1})^T \in \rr^{M_1}$ and $\boldsymbol{\zeta}^{(+)}:=(\zeta_1^{(+)},\ldots,\zeta^{(+)}_{M_2})^T \in \rr^{M_2}$.  Then using  (\ref{se1}) and (\ref{se2}),  we write the following upper bounds  for the $\ell_2$-norms of the vectors  $\boldsymbol{\rho}^{(-)}$ and $\boldsymbol{\rho}^{(+)}$
\begin{align}
\|\boldsymbol{\rho}^{(-)}\|_2  \leq \|\boldsymbol{\zeta}^{(-)}\|_2 \,   \kappa^{2} \left( \mathbf{V}_{M_1} (\boldsymbol{z}^{(-)}) \right),  \quad
\|\boldsymbol{\rho}^{(+)}\|_2  \leq \|\boldsymbol{\zeta}^{(+)}\|_2 \,   \kappa^{2} \left( \mathbf{V}_{M_2} (\boldsymbol{z}^{(+)}) \right) . \label{senb}
\end{align}
Applying ideas similar to those in Subsection \ref{unsen},  we write also the upper bounds for the $\ell_1$-norms of the sensitivity vectors $\boldsymbol{\rho}^{(-)}$ and $\boldsymbol{\rho}^{(+)}$
\begin{align}
\| \boldsymbol{\rho}^{(-)}\|_1 & \leq \zeta^{(-)}_{\max} \,   \|\mathbf{V}_{M_1} (\boldsymbol{z}^{(-)})  \|_2^{2} \,   \|\mathbf{V}^{-1}_{M_1}(\boldsymbol{z}^{(-)})    \|_{\mathrm{F}}^{2}, \label{senb1}  \\
\| \boldsymbol{\rho}^{(+)}\|_1 &
 \leq \zeta^{(+)}_{\max} \,   \|\mathbf{V}_{M_2} (\boldsymbol{z}^{(+)})  \|_2^{2} \,   \|\mathbf{V}^{-1}_{M_2}(\boldsymbol{z}^{(+)})   \|_{\mathrm{F}}^{2} , \label{senb2}
\end{align}
where  $\zeta_{\max}^{(-)}:=\max_{j=1,\ldots,M_1} \zeta_j^{(-)}$ and $\zeta_{\max}^{(+)}:=\max_{j=1,\ldots,M_2} \zeta_j^{(+)}$.  Naturally, arises the question of how sharp are the inequalities (\ref{senb})--(\ref{senb2}). This question is discussed  in Example \ref{ex2}.

Let us now consider  the case of structured perturbations.  Define further  matrices  
$$
\mathbf{\mathcal{S}}^{(-)}_{M_1,2M_1}:=\big(S_{n,m}^{(-)} \big)_{n,m=1}^{M_1,2M_1} \in \cc^{M_1\times 2M_1},   \quad  \mathbf{\mathcal{S}}^{(+)}_{M_2,2M_2}:=\big(S_{n,m}^{(+)}\big)_{n,m=1}^{M_2,2M_2} \in \cc^{M_2\times 2M_2},
$$
 and vectors
 \begin{align*}
 (\boldsymbol{\mathcal{S}}^{(n)}_{2M_1})^{(-)}&:=\left((S^{(n)}_1)^{(-)}, \ldots, (S^{(n)}_{2M_1})^{(-)} \right)^{T} \in \cc^{2M_1} \quad \text{for} \quad n=1,\ldots,M_1,  \\
 (\boldsymbol{\mathcal{S}}^{(n)}_{2M_2})^{(+)}&:=\left((S^{(n)}_1)^{(+)}, \ldots, (S^{(n)}_{2M_2})^{(+)}\right)^{T} \in \cc^{2M_2} \quad \text{for} \quad n=1,\ldots,M_2,  
 \end{align*}
  by their components
 \begin{align}
S_{n,m}^{(-)}=\left(S^{(n)}_{m} \right)^{(-)}: =& \frac{1}{\gamma_{n}^{(-)}} \bigg( \,  \sum\limits_{k=\max\{1,m-M_1\}}^{\min\{M_1,m-1 \}} \bigg( \left( \mathbf{V}^{T}_{M_1} (\boldsymbol{z}^{(-)}) \right)^{-1}\bigg)_{k,n} \,  \bigg( \left( \mathbf{V}^{T}_{M_1} (\boldsymbol{z}^{(-)})\right)^{-1} \bigg)_{m-k,n} \notag \\
& - z_{n}^{(-)} \,  \sum\limits_{k=\max\{1,m-M_1+1\}}^{\min\{M_1,m \}} \bigg( \left( \mathbf{V}^{T}_{M_1} (\boldsymbol{z}^{(-)}) \right)^{-1} \bigg)_{k,n} \, \bigg(\left( \mathbf{V}^{T}_{M_1} (\boldsymbol{z}^{(-)}) \right)^{-1} \bigg)_{m-k+1,n} \bigg)\label{sent1}
\end{align}
and
 \begin{align}
S_{n,m}^{(+)}=\left(S^{(n)}_{m} \right)^{(+)} :=& -\frac{1}{ \gamma_{n}^{(+)} \, z_{n}^{-2}} \bigg( \,  \sum\limits_{k=\max\{1,m-M_2\}}^{\min\{M_2,m-1 \}} \bigg( \left( \mathbf{V}^{T}_{M_2} (\boldsymbol{z}^{(+)}) \right)^{-1} \bigg)_{k,n} \,  \bigg( \left( \mathbf{V}^{T}_{M_2} (\boldsymbol{z}^{(+)})\right)^{-1} \bigg)_{m-k,n} \notag \\
& - \frac{1}{z_{n}^{(+)} }  \,  \sum\limits_{k=\max\{1,m-M_2+1\}}^{\min\{M_2,m \}} \bigg( \left( \mathbf{V}^{T}_{M_2} (\boldsymbol{z}^{(+)}) \right)^{-1} \bigg)_{k,n} \, \bigg( \left( \mathbf{V}^{T}_{M_2} (\boldsymbol{z}^{(+)}) \right)^{-1} \bigg)_{m-k+1,n} \bigg). \label{sent2}
\end{align}
Formulas (\ref{sent1}) and (\ref{sent2}) are obtained using  ideas similar to those presented in Subsection \ref{stsen}.  In (\ref{sent2})  we take into account (\ref{gam1}).

Now we define the  sensitivity $\eta_{n,m}^{(-)}$ of the pole $z^{(-)}_{n}$ located inside the unit circle with respect to the Fourier coefficient $\hat{r}_{-m}$ by $\eta_{n,m}^{(-)}:=|S_{n,m}^{(-)}|$  for $n=1,\ldots,M_1$ and $m=1,\ldots,2M_1$.  Similarly,  by $\eta_{n,m}^{(+)}:=|S_{n,m}^{(+)}|$ we define the sensitivity of the pole $z^{(+)}_{n}$ located outside the unit circle with respect to the Fourier coefficient $\hat{r}_{m}$ for $n=1,\ldots,M_2$ and $m=1,\ldots,2M_2$. By 
\begin{equation}\label{strvect1}
   (\boldsymbol{\eta}^{(n)}_{2M_1})^{(-)}:=(\eta_{n,1}^{(-)},\ldots, \eta_{n,2M_1}^{(-)})^T \text{ and }  (\boldsymbol{\eta}^{(n)}_{2M_2})^{(+)}:=(\eta_{n,1}^{(+)},\ldots, \eta_{n,2M_2}^{(+)})^T 
\end{equation}
 we denote the corresponding sensitivity vectors. Note that Fourier coefficients $\hat{r}_{-m}$ and $\hat{r}_{m}$ are $m$th measurements in this case.  We consider also two vectors 
$$
\boldsymbol{\eta}^{(-)}:=(\eta_1^{(-)},\ldots, \eta_{M_1}^{(-)})^{T} \in \mathbb{R}^{M_1},  \quad  \boldsymbol{\eta}^{(+)}:=(\eta_1^{(+)},\ldots, \eta_{M_2}^{(+)})^{T} \in \mathbb{R}^{M_2},
$$ 
where their components
\begin{align}
    \eta_n^{(-)} &:=\|(\boldsymbol{\mathcal{S}}^{(n)}_{2M_1})^{(-)}\|_2=\sqrt{(\eta^{(-)}_{n,1})^{2}+\ldots+(\eta^{(-)}_{n,2M_1})^{2}},  \quad n=1,\ldots,M_1,  \label{stp1}  \\
    \eta_n^{(+)} &:=\|(\boldsymbol{\mathcal{S}}^{(n)}_{2M_2})^{(+)}\|_2=\sqrt{(\eta^{(+)}_{n,1})^{2}+\ldots+(\eta^{(+)}_{n,2M_2})^{2}},  \quad n=1,\ldots,M_2, \label{stp2}
\end{align}
are sensitivities of the poles $z_{n}^{(-)}$ and $z_{n}^{(+)}$, respectively, with respect to structured perturbations.

\begin{table}[h!]
\centering
 \caption{\small   Reconstruction errors $|z_j-\tilde{z}_j|$ for poles $z_j$    for Example \ref{ex3} }
  \label{tab_nst1ex1}
\begin{tabular}{ |p{1.7cm}||p{2.2cm}|p{2.2cm}||  p{2.2cm}|p{2.2cm}|}
 \hline
$\sigma=10^{-3}$   & \multicolumn{2}{|c||}{$N =2$} & \multicolumn{2}{|c|}{$N =10$}\\
 \cline{2-5}
 & $z_1^{(-)}=-0.1$ & $z_1^{(+)}=-2.1$  &  $z_1^{(-)}=-0.1$   &$z_1^{(+)}=-2.1$  \\[2mm]
 \hline
min & $2.54$e--$05$    &$2.06$e--$04$ & $3.49$e--$06$ & $5.56$e--$05$  \\[2mm]
max  & $2.85$e--$04$    &$3.34$e--$03$ & $2.7$e--$04$ & $3.15$e--$03$   \\[2mm]
average & $1.34$e--$04$    &$2.03$e--$03$ &  $1.29$e--$04$ & $1.06$e--$03$   \\
 \hline
 \end{tabular}

\begin{tabular}{ |p{1.7cm}||p{2.2cm}|p{2.2cm}|| p{2.2cm}|p{2.2cm}|}
 \hline
 $\sigma=10^{-1}$   & \multicolumn{2}{|c||}{$N =10$} & \multicolumn{2}{|c|}{$N =30$}\\
    \cline{2-5}
 & $z_1^{(-)}=-0.1$ & $z_1^{(+)}=-2.1$  &  $z_1^{(-)}=-0.1$   &$z_1^{(+)}=-2.1$  \\[2mm]
 \hline
min & $3.2$e--$03$    &$4.09$e--$02$ & $7.19$e--$03$ & $4.15$e--$02$  \\[2mm]
max  & $2.28$e--$02$    &$3.65$e--$01$ & $2.72$e--$02$ & $1.7$e--$01$   \\[2mm]
average & $1.09$e--$02$    &$1.57$e--$01$ &  $1.36$e--$02$ & $1.05$e--$01$   \\
 \hline
 \end{tabular}
  
 \end{table}

\section{Numerical experiments}
In this section,  we present several numerical experiments.
The algorithms are implemented in \textsc{Matlab} and use IEEE standard floating point arithmetic with double precision. 
By $\tilde{z}_j$ and $\tilde{\gamma}_j$ we denote reconstructed   poles  and  coefficients, respectively, and by the following formulas
$$
e(\boldsymbol{z}) \coloneqq \max\limits_{j=1,\ldots,M}|z_j-\tilde{z}_j|,
\quad e(\gamra) \coloneqq \max\limits_{j=1,\ldots,M}|\gamma_j-\tilde{\gamma}_j|,
$$
we define the corresponding reconstruction errors.
  Note that in the case when we deal with noisy data,  we assume that a number  of poles $M$,  as well as numbers of poles located inside and outside the unit circle  $M_1$ and $M_2$, respectively, are known.

   \begin{table}[h!]
\centering
 \caption{\small Reconstruction errors $e(\gamra)$   for Example \ref{ex3}  }
  \label{tab_coef}
\begin{tabular}{ |p{1.7cm}||p{1.55cm}|p{1.55cm}||p{1.7cm}||p{1.55cm}|p{1.55cm}| }
   \hline
$\sigma=10^{-3}$ & $N=2$ & $N=10$ & $\sigma=10^{-1}$ &$N =10$  & $N=30$ \\[2mm]
 \hline
min & $1.78$e--$04$     &$1.17$e--$04$ & min &  $4.03$e--$02$    &$2.77$e--$02$  \\[2mm]
max & $2.11$e--$03$    &$2.46$e--$03$ & max &  $3.04$e--$01$    &$1.44$e--$01$  \\[2mm]
average & $1.44$e--$03$    &$8.48$e--$04$ & average &  $1.11$e--$01$    &$9.55$e--$02$  \\
 \hline
 \end{tabular}
\end{table}

  In Examples \ref{ex3}-\ref{ex2}, we demonstrate recovery results for three different rational functions and analyze them employing sensitivity analysis as was discussed  in Section \ref{method}. 
  We utilize $2N$ noisy Fourier coefficients  $\hat{r}_{-k}(1 + \varepsilon_{k})$, $k=1, \ldots, 2N$,  to reconstruct poles $z_n^{(-)}$,   and $2N$ noisy Fourier coefficients $\hat{r}_{k}(1 + \varepsilon_{k})$, $k=1, \ldots, 2N$,  to reconstruct poles $z_n^{(+)}$.  
By $(\varepsilon_{k})_{k=1}^{2N}$ we denote the random noise that is generated in \textsc{Matlab} by 
\texttt{$\sigma$*randn(2N)},  i.e.  we consider a multiplicative  noise drawn from the standard normal distribution.  We set $L=N$ and perform computations for different values of $N$ and  $\sigma$.  For $\sigma>0$ we perform 10 itterations.
To compute sensitivities $\rho_n^{(-)}$,  $\rho_n^{(+)}$ and ,$\eta_n^{(-)}$,   $\eta_n^{(+)}$,  we use formulas  (\ref{senk}),  (\ref{senk2}) and (\ref{stp1}),  (\ref{stp2}),  respectively.   Upper bounds for sensitivities $\rho_n^{(-)}$ and $\rho_n^{(+)}$  are computed via (\ref{se1}) and (\ref{se2}), respectively.

\begin{table}[h!]
\centering
 \caption{\small Reconstruction errors $|z_j-\tilde{z}_j|$ for poles $z_j$ for Example \ref{ex1} for $\sigma=10^{-3}$ }
  \label{tab_nst52}
\begin{tabular}{ |p{1.3cm}||p{2.2cm}|p{2.2cm}| p{2.2cm}|p{2.2cm}|}
   \hline
$N =8$ & $z_1^{(-)}=0.9$ & $z_2^{(-)}=-0.9$  & $z_3^{(-)}=0.9\mathrm{i}$  &$z_4^{(-)}=-0.9\mathrm{i}$  \\[2mm]
 \hline
min & $7.08$e--$05$    &$4.31$e--$06$ & $2.79$e--$06$ & $9.06$e--$06$  \\[2mm]
max  & $6.74$e--$04$    &$3.17$e--$04$ & $1.65$e--$04$ & $1.45$e--$04$   \\[2mm]
average & $3.18$e--$04$    &$1.48$e--$04$ & $7.87$e--$05$ & $5.57$e--$05$   \\[2mm]
 \hline
 \end{tabular}
  \qquad 
  \begin{tabular}{ |p{1.3cm}||p{2.2cm}|p{2.2cm}| p{2.2cm}|p{2.2cm}|}
   \hline
$N =8$ & $z_1^{(+)}=1.1$ & $z_2^{(+)}=-1.1$  & $z_3^{(+)}=1.1\mathrm{i}$  &$z_4^{(+)}=-1.1\mathrm{i}$  \\[2mm]
 \hline
min & $1.71$e--$05$    &$1.13$e--$06$ & $1.31$e--$05$ & $8.74$e--$06$  \\[2mm]
max  & $2.86$e--$04$    &$2.7$e--$04$ & $1.89$e--$04$ & $1.9$e--$04$   \\[2mm]
average & $1.44$e--$04$    &$1.19$e--$04$ & $1.04$e--$04$ & $9.43$e--$05$   \\[2mm]
 \hline
 \end{tabular}
 \end{table}

\begin{table}[h!]
\centering
 \caption{\small Reconstruction errors $|z_j-\tilde{z}_j|$ for poles $z_j$ for Example \ref{ex1} for $\sigma=10^{-1}$ }
  \label{tab_nst521}
\begin{tabular}{ |p{1.3cm}||p{2.2cm}|p{2.2cm}| p{2.2cm}|p{2.2cm}|}
   \hline
$N =8$ & $z_1^{(-)}=0.9$ & $z_2^{(-)}=-0.9$  & $z_3^{(-)}=0.9\mathrm{i}$  &$z_4^{(-)}=-0.9\mathrm{i}$  \\[2mm]
 \hline
min & $6.41$e--$03$    &$2.65$e--$03$ & $5.14$e--$04$ & $1.38$e--$03$  \\[2mm]
max  & $2.02$e--$01$    &$2.37$e--$02$ & $2.89$e--$02$ & $7.51$e--$03$   \\[2mm]
average & $3.99$e--$02$    &$1.13$e--$02$ & $1.06$e--$02$ & $1.82$e--$02$   \\[2mm]
 \hline
 \end{tabular}
  \qquad 
  \begin{tabular}{ |p{1.3cm}||p{2.2cm}|p{2.2cm}| p{2.2cm}|p{2.2cm}|}
   \hline
$N =8$ & $z_1^{(+)}=1.1$ & $z_2^{(+)}=-1.1$  & $z_3^{(+)}=1.1\mathrm{i}$  &$z_4^{(+)}=-1.1\mathrm{i}$  \\[2mm]
 \hline
min & $8.17$e--$04$    &$8.71$e--$04$ & $1.93$e--$03$ & $2.83$e--$03$  \\[2mm]
max  & $2.95$e--$02$    &$2.61$e--$02$ & $2.48$e--$02$ & $2.41$e--$02$   \\[2mm]
average & $1.56$e--$02$  &  $1.4$e--$02$ & $1.24$e--$02$ & $1.15$e--$02$   \\[2mm]
 \hline
 \end{tabular}
 \end{table}

  \begin{table}[h!]
\centering
 \caption{\small Reconstruction errors $e(\gamra)$   for Example \ref{ex1}  }
  \label{tab_coef52}
\begin{tabular}{ |p{1.7cm}||p{1.8cm}|p{1.8cm}|}
   \hline
$N=8$ & $\sigma=10^{-3}$ & $\sigma=10^{-1}$  \\[2mm]
 \hline
min & $4.2$e--$03$     &$4.25$e--$01$   \\[2mm]
max & $1.73$e--$02$    &$1.22$e--$00$   \\[2mm]
average & $9.7$e--$03$    &$6.35$e--$01$   \\
 \hline
 \end{tabular}
\end{table}

\begin{example}\label{ex3}
 Let a rational function $r(z)$ as in (\ref{rat}) be given by parameters $M=2$ and
$$
\boldsymbol{z}=( -0.1,  -2.1)^{T},  \quad  \boldsymbol{\gamma}=(0.5,  0.5)^{T}
$$
(see \cite{EI22,  ZGA21}).   
 When $\sigma=0$,  i.e.  in case of recovery from exact data,  we get the following reconstruction errors for $N=2$
$$
e(\boldsymbol{z}) =4.44 \cdot 10^{-16},  \quad e(\gamra)  = 1.11 \cdot 10^{-16}.
$$
  Since we are  interested to investigate which pole is the most sensitive to perturbations,   we compute reconstruction errors for poles $z_1^{(-)}=-0.1$ and $z_1^{(+)}=-2.1$ independently and present the results in Table \ref{tab_nst1ex1}.  Results for the recovery of the coefficient vector $\gamra$ are  shown in Table \ref{tab_coef}.

\begin{table}[h!]
\centering
\caption{\small {Sensitivities for Example \ref{ex1}.} }
\label{tb1} 
\begin{tabular}{ |p{0.2cm}||p{2cm}|  p{1.5cm}|  p{1.5cm}|  p{2.2cm}|  }
 \hline
$j$ &\ \ \ \ poles $z_j^{(-)}$ &\ \ \ \  $\eta_j^{(-)}$  &\ \ \ \  $\rho_j^{(-)}$   &\ \ \ \   bound $\rho_j^{(-)}$ \\[2mm]
 \hline
1 &   $0.9$ &   0.243
    & 7.748  &  
  13.548
  \\[2mm]
2 &   $-0.9$ &    0.121   &3.874  &   6.774  \\[2mm]
3 &   $0.9 \mathrm{i}$ & 0.081  &  2.582 &   4.516 \\[2mm]
4 &    $-0.9 \mathrm{i}$ & 0.06  &1.937  &   3.387  \\[2mm]
 \hline
 \end{tabular}
 \qquad 
\begin{tabular}{ |p{0.2cm}||p{2cm}|   p{1.5cm}| p{1.5cm}|  p{2.2cm}|  }
 \hline
$j$ &\ \ \ \ poles $z_j^{(+)}$ &\ \ \ \  $\eta_j^{(+)}$    &\ \ \ \  $\rho_j^{(+)}$   &\ \ \ \   bound $\rho_j^{(+)}$ \\[2mm]
 \hline
1 &  $1.1$  &    0.056
     &3.533
  &   5.695
    \\[2mm]
2 &   $-1.1$ &    0.047  & 2.944  &   4.745  \\[2mm]
3 &   $1.1 \mathrm{i}$ &  0.04  & 2.523  &   4.067  \\[2mm]
4 &    $-1.1 \mathrm{i}$ &   0.035    & 2.208  &   3.559  \\[2mm]
 \hline
 \end{tabular}
\end{table}
 
 Analyzing numerical  results from Table \ref{tab_nst1ex1},  we conclude that for  $\sigma=10^{-1}$ our method still reconstructs the locations of the poles quite good.  On the other hand,  employing more noisy Fourier coefficients also gives us better recovery results.   We notice 
  that a pole $z_1^{(-)}=-0.1$ is less sensitive to perturbations than a pole $z_1^{(+)}=-2.1$.  This numerical observation is also confirmed with the sensitivities values.  For a pole $z_1^{(-)}=-0.1$,  we have $\rho_1^{(-)}=2\cdot 10^{-1}$ and  $\eta_1^{(-)}=1.97$,  and for a pole $z_1^{(+)}=-2.1$,  we have $\rho_1^{(+)}=9.52\cdot 10^{-1}$ and  $\eta_1^{(+)}=6.65$.  Both unstructured and structured sensitivities are larger for a pole $z_1^{(+)}=-2.1$  located outside the unit circle.  We compute also vectors (\ref{strvect1})
 $$
 (\boldsymbol{\eta}^{(1)}_{2})^{(-)}=(2.0 \cdot 10^{-1},  2.0)^{T} \text{  and  }  (\boldsymbol{\eta}^{(1)}_{2})^{(+)}=(4.2,  8.82)^{T},
 $$
which provide us with the additional information that a pole  $z_1^{(-)}=-0.1$ is the most sensitive to the perturbations of  the Fourier coefficient $\hat{r}_{-2}$,  and a pole $z_1^{(+)}=-2.1$ is the most sensitive  to the perturbations of  the Fourier coefficient $\hat{r}_{2}$.
\end{example}

\begin{table}[h!]
\centering
 \caption{\small Reconstruction errors $|z_j-\tilde{z}_j|$ for poles $z_j$ for Example \ref{ex2} for $N=4$ }
  \label{tab_nst531}
\begin{tabular}{ |p{1.8cm}||p{2.2cm}|p{2.2cm}| p{2.2cm}|p{2.2cm}|}
   \hline
$\sigma=10^{-7}$  & $z_1^{(-)}=0.2$ & $z_2^{(-)}=0.5$  & $z_1^{(+)}=2$  &$z_2^{(+)}=50$  \\[2mm]
 \hline
min & $1.58$e--$08$    &$2.9$e--$08$ & $6.07$e--$09$ & $1.75$e--$03$  \\[2mm]
max  & $5.1$e--$07$    &$2.49$e--$07$ & $2.7$e--$07$ & $3.21$e--$01$   \\[2mm]
average & $2.48$e--$07$    &$1.3$e--$07$ & $1.45$e--$07$ & $1.34$e--$01$   \\[2mm]
 \hline
   \hline
$\sigma=10^{-5}$ & $z_1^{(-)}=0.2$ & $z_2^{(-)}=0.5$  & $z_1^{(+)}=2$  &$z_2^{(+)}=50$  \\[2mm]
 \hline
min & $4.24$e--$06$    &$1.18$e--$06$ & $2.57$e--$06$ & $3.16$e--$00$  \\[2mm]
max  & $6.45$e--$05$    &$2.77$e--$05$ & $2.92$e--$05$ & $9.46$e+$01$   \\[2mm]
average & $2.31$e--$05$    &$1.13$e--$05$ & $1.5$e--$05$ & $2.02$e+$01$   \\[2mm]
 \hline
 \end{tabular}
  \end{table}

\begin{example}\label{ex1}

We consider the case when poles are located very close to the unit circle.  Let a rational function $r(z)$ as in (\ref{rat}) of type $(7,8)$,  i.e with $M=8$,  be given by vectors of poles and coefficients
$$
\boldsymbol{z}=(0.9, -0.9,  0.9 \mathrm{i},  -0.9 \mathrm{i},  1.1,  -1.1,  1.1 \mathrm{i},  -1.1 \mathrm{i})^{T},  \quad  \boldsymbol{\gamma}=(1,2,3,4,5,6,7,8)^{T},
$$
respectively. First,  we consider  recovery from the exact data,  i.e for $\sigma=0$.  For  $N=8$, we obtain the following reconstruction errors
$$
e(\boldsymbol{z})=1.45\cdot 10^{-15},  \quad e( \boldsymbol{\gamma})=5.48 \cdot 10^{-14}.
$$

To reconstruct from noisy data, we consider $\sigma=10^{-3}$ and $\sigma=10^{-1}$. The recovery results for each pole $z_j$ are presented in Tables \ref{tab_nst52} and \ref{tab_nst521}, as well as for the coefficient vector $\gamra$ in Table \ref{tab_coef52}.  Based on these numerical results,  we conclude that a pole $z_1^{(-)}=0.9$ is the most sensitive to perturbations, which is corroborated  by the sensitivities shown in Table \ref{tb1}.

\end{example}

 \begin{table}[h!]
\centering
\caption{\small {Sensitivities for Example \ref{ex2}.} }
\label{tb53} 
\begin{tabular}{ |p{0.2cm}||p{2cm}|  p{2cm}|  p{2cm}|  p{2.2cm}|  }
 \hline
$j$ &\ \ \ \ poles $z_j^{(-)}$ &\ \ \ \  $\eta_j^{(-)}$  &\ \ \ \  $\rho_j^{(-)}$   &\ \ \ \   bound $\rho_j^{(-)}$ \\[2mm]
 \hline
1 &   $0.2$ &   $1.807 \cdot 10^{1}$
    & $1.766 \cdot 10^{1}$  &  
  $3.937 \cdot 10^{1}$
  \\[2mm]
2 &   $0.5$ &    $1.518 \cdot 10^{1}$   &$2.249 \cdot 10^{1}$ &  $5.625 \cdot 10^{1}$  \\[2mm]
 \hline
 \end{tabular}
 \qquad 
\begin{tabular}{ |p{0.2cm}||p{2cm}|   p{2cm}| p{2cm}|  p{2.2cm}|  }
 \hline
$j$ &\ \ \ \ poles $z_j^{(+)}$ &\ \ \ \  $\eta_j^{(+)}$    &\ \ \ \  $\rho_j^{(+)}$   &\ \ \ \   bound $\rho_j^{(+)}$ \\[2mm]
 \hline
1 &  $2$  &    $1.973 \cdot 10^{1}$
     & 5.43
  &   $1.993 \cdot 10^{1}$
    \\[2mm]
2 &   $50 $ &    $1.577 \cdot 10^{4}$  & $2.204 \cdot 10^{3}$  &   $6.477 \cdot 10^{3}$  \\[2mm]
 \hline
 \end{tabular}
\end{table}

\begin{example}\label{ex2}

Finally,  let a rational function $r(z)$ as  in (\ref{rat}) of type $(3,4)$,  i.e with $M=4$,  be given by vectors of poles and coefficients
$$
\boldsymbol{z}=(0.2,  0.5, 2,  50)^{T},  \quad  \boldsymbol{\gamma}=(1,1,1, 1)^{T},
$$
respectively.   For  $N=4$ and $\sigma=0$,  we get the following reconstruction errors
$$
e(\boldsymbol{z})=1.42\cdot 10^{-13},  \quad e( \boldsymbol{\gamma})=9.27 \cdot 10^{-15}.
$$

Then we consider recovery from noisy data with $\sigma=10^{-7}$ and $\sigma=10^{-5}$, and present the corresponding results in Table \ref{tab_nst531}.   Analyzing these numerical results,  we can conclude that a pole $z_2^{(+)}=50$ is the most sensitive to perturbations.  For $\sigma=10^{-5}$ the method does not reconstruct location of this pole correctly.   This numerical observation is also predicted by the sensitivities values shown in Table \ref{tb53}.  Both unstructured $\rho_2^{(+)}=2.204 \cdot 10^{3}$ and structured $\eta_2^{(+)}=1.577 \cdot 10^{4}$ sensitivities are sufficiently larger than sensitivities of other poles.   Note,  that the locations of the poles $z_1^{(-)}=0.2$,  $z_2^{(-)}=0.5$ and  $z_1^{(+)}=2$  are reconstructed  even for $\sigma=10^{-1}$.  We present results regarding recovery of  the coefficient vector $\gamra$ in Table \ref{tab_coef53}.

  \begin{table}[h!]
\centering
 \caption{\small Reconstruction errors $e(\gamra)$   for Example  \ref{ex2} }
  \label{tab_coef53}
\begin{tabular}{ |p{1.7cm}||p{1.8cm}|p{1.8cm}|}
   \hline
$N=4$ & $\sigma=10^{-7}$ & $\sigma=10^{-5}$  \\[2mm]
 \hline
min & $7.38$e--$04$     &$3.79$e--$02$   \\[2mm]
max & $7.38$e--$03$    &$6.95$e--$01$   \\[2mm]
average & $4$e--$03$    &$4.11$e--$01$   \\
 \hline
 \end{tabular}
\end{table}

We compute also vectors (\ref{strvect1}) for poles located inside the unit circle
$$
 (\boldsymbol{\eta}^{(1)}_{4})^{(-)}=(0.555,  4.999, 13.333, 11.111)^{T},   \quad  (\boldsymbol{\eta}^{(2)}_{4})^{(-)}=(0.222,  2.666, 9.999, 11.111)^{T},  
$$
and outside the unit circle
\begin{align*}
     (\boldsymbol{\eta}^{(1)}_{4})^{(+)}&=(3.472\cdot 10^{-3},  0.354, 9.375, 17.36)^{T}, \\
     (\boldsymbol{\eta}^{(2)}_{4})^{(+)}&=(54.25,  2.929 \cdot 10^{3}, 1.106 \cdot 10^{4}, 1.085 \cdot 10^{4})^{T}.
\end{align*}
This implies that the poles $z_1^{(-)}=0.2$ and $z_2^{(-)}=0.5$ are the most sensitive to perturbations of the Fourier coefficients $\hat{r}_{-3}$ and $\hat{r}_{-4}$, respectively, while the poles $z_1^{(+)}=2$ and $z_2^{(+)}=50$ are the most sensitive to the perturbations of the Fourier coefficients $\hat{r}_{4}$ and $\hat{r}_{3}$, respectively.

\begin{table}[h!]
\centering
 \caption{\small  $\ell_p$-norms  of  $\boldsymbol{\rho}^{(-)}$ and $\boldsymbol{\rho}^{(+)}$ and their upper bounds for Example \ref{ex2} for $p=1,2$  }
 \label{tab_norm1}
\begin{tabular}{ |p{1.2cm}||p{1.5cm}|p{2.8cm}| p{2cm}|p{2.8cm}| }
   \hline
 & $\|\boldsymbol{\rho}^{(-)}\|_2$  & bound ($\|\boldsymbol{\rho}^{(-)}\|_2$)  & $\|\boldsymbol{\rho}^{(+)}\|_2$ & bound ($\|\boldsymbol{\rho}^{(+)}\|_2$) \\[2mm]
 \hline
$N=4$ & $2.86$e+$01$    &$6.866$e+$01$ & $2.204$e+$03$ &   $6.477$e+$03$    \\[2mm]
 \hline
   \hline
 & $\|\boldsymbol{\rho}^{(-)}\|_1$  & bound ($\|\boldsymbol{\rho}^{(-)}\|_1$)  & $\|\boldsymbol{\rho}^{(+)}\|_1$ & bound ($\|\boldsymbol{\rho}^{(+)}\|_1$) \\[2mm]
 \hline
$N=4$ & $4.01$e+$01$    &$5.725$e+$01$ & $2.209$e+$03$ &   $6.802$e+$03$    \\[2mm]
 \hline
 \end{tabular}
\end{table}

In Subsection \ref{sub42}, we derived  the upper bounds  (\ref{senb})--(\ref{senb2})  for the $\ell_2$- and the $\ell_1$-norms  of vectors $\boldsymbol{\rho}^{(-)}$ and $\boldsymbol{\rho}^{(+)}$. In order to understand how sharp these estimates are, we compute the exact values for the corresponding norms and their upper bounds as in  (\ref{senb})--(\ref{senb2}) and present the obtained numerical results in Table \ref{tab_norm1}. Comparing these results, we can conclude that the estimates (\ref{senb})--(\ref{senb2}) are quite sharp.

 \begin{table}[h!]
\centering
 \caption{\small 
 Reconstruction results for Example  \ref{ex3} }
  \label{tab_16}
\begin{tabular}{ |p{1.7cm}||p{1.5cm}|p{1.5cm}||p{1.7cm}||p{1.5cm}|p{1.5cm}| }
 \hline
    \multicolumn{3}{|c||}{Method via Hankel pencil} & \multicolumn{3}{|c|}{AAA method}\\
   \hline
$\sigma=10^{-7}$ & $e(\boldsymbol{z})$ & $e(\gamra) $ & $\sigma=10^{-7}$ & $e(\boldsymbol{z})$  & $e(\gamra) $ \\[2mm]
 \hline
min & $1.52$e--$08$     &$1.23$e--$08$ & min &  $3.54$e--$08$    &$4.55$e--$08$  \\[2mm]
max & $1.19$e--$07$    &$5.51$e--$08$ & max &  $6.62$e--$07$    &$4.84$e--$07$  \\[2mm]
average & $7.51$e--$08$    &$2.89$e--$08$ & average &  $3.83$e--$07$    &$2.82$e--$07$  \\
 \hline
 \end{tabular}
 \begin{tabular}{ |p{1.7cm}||p{1.5cm}|p{1.5cm}||p{1.7cm}||p{1.5cm}|p{1.5cm}| }
 \hline
    \multicolumn{3}{|c||}{Method via Hankel pencil} & \multicolumn{3}{|c|}{AAA method}\\
   \hline
$\sigma=10^{-3}$ & $e(\boldsymbol{z})$ & $e(\gamra) $ & $\sigma=10^{-3}$ & $e(\boldsymbol{z})$  & $e(\gamra) $ \\[2mm]
 \hline
min & $4.34$e--$03$     &$1.06$e--$03$ & min &  $2.37$e--$02$    &$2.27$e--$02$  \\[2mm]
max & $1.23$e--$02$    &$4.24$e--$03$ & max &  $2.05$e+$00$    &$7.56$e--$02$  \\[2mm]
average & $7.64$e--$03$    &$2.71$e--$03$ & average &  $4.52$e--$01$    &$4.55$e--$02$  \\
 \hline
 \end{tabular}
\end{table}

\begin{table}[h!]
\centering
 \caption{\small Reconstruction results   for  Example \ref{ex1} }
  \label{tab_17}
\begin{tabular}{ |p{1.7cm}||p{1.5cm}|p{1.5cm}||p{1.7cm}||p{1.5cm}|p{1.5cm}| }
 \hline
    \multicolumn{3}{|c||}{Method via Hankel pencil} & \multicolumn{3}{|c|}{AAA method}\\
   \hline
$\sigma=10^{-4}$ & $e(\boldsymbol{z})$ & $e(\gamra) $ & $\sigma=10^{-4}$ & $e(\boldsymbol{z})$  & $e(\gamra) $ \\[2mm]
 \hline
min & $1.96$e--$05$     &$4.79$e--$06$ & min &  $7.03$e--$05$    &$7.18$e--$05$  \\[2mm]
max & $1.37$e--$04$    &$6.22$e--$05$ & max &  $4.49$e--$04$    &$3.43$e--$04$  \\[2mm]
average & $5.89$e--$05$    &$2.34$e--$05$ & average &  $2.98$e--$04$    &$2.25$e--$04$  \\
 \hline
 \end{tabular}
 \begin{tabular}{ |p{1.7cm}||p{1.5cm}|p{1.5cm}||p{1.7cm}||p{1.5cm}|p{1.5cm}| }
 \hline
    \multicolumn{3}{|c||}{Method via Hankel pencil} & \multicolumn{3}{|c|}{AAA method}\\
   \hline
$\sigma=10^{-2}$ & $e(\boldsymbol{z})$ & $e(\gamra) $ & $\sigma=10^{-2}$ & $e(\boldsymbol{z})$  & $e(\gamra) $ \\[2mm]
 \hline
min & $2.2$e--$03$     &$2.42$e--$03$ & min &  $1.54$e--$02$    &$9.73$e--$03$  \\[2mm]
max & $1.37$e--$02$    &$4.86$e--$02$ & max &  $6.69$e--$02$    &$3.54$e--$02$  \\[2mm]
average & $6.93$e--$03$    &$1.06$e--$02$ & average &  $4.42$e--$02$    &$5.75$e--$02$  \\
 \hline
 \end{tabular}
\end{table}

\end{example}

  In Example \ref{ex4}, we discuss the comparison of the performance of our approach with the AAA method for rational approximation.  We employ $4N$ noisy function values  $r(t_k)(1 + \varepsilon_{k})$ with $t_k=\mathrm{e}^{\frac{2\pi \mathrm{i} k}{4N}}$ for $k=0, \ldots, 4N-1$, as input information for both methods. By $(\varepsilon_{k})_{k=0}^{4N-1}$ as above we denote the random noise that is generated in \textsc{Matlab} by 
\texttt{$\sigma$*randn(4N)}.  The corresponding modification of our method is discussed in Remark \ref{rem41}. We refer the interested reader to \cite{AAA} regarding description of the AAA method. Note only that
 AAA takes the form of iteration with respect to a type of a rational approximant that is represented    at each step in the barycentric form. In  case of reconstruction of rational function of type $(M-1,M)$ using exact data,  AAA  terminates after $M+1$ iteration steps and returns  $r$ as
$$
r(t)= \sum\limits_{s=1}^{M+1}  \frac{r(t_{k_s})w_s}{t-t_{k_s}} \bigg/ \sum\limits_{s=1}^{M+1} \frac{w_s}{t-t_{k_s}},
$$
 where points $t_{k_s}$ are chosen by some greedy algorithm and weights $w_s$ are computed by
  solving a linear least-squares problem over the subset of sample points that have not been selected as support points,  $ \{t_{0},\ldots,t_{4N-1} \}\setminus \{t_{k_1},\ldots,t_{k_{M+1}} \}$ (see also \cite{DPP21, DPR23, DP21} for details). The  AAA method has the overall computational cost of $\mathcal{O}(N M^3)$ operations.

\begin{table}[h!]
\centering
 \caption{\small Reconstruction results   for  Example \ref{ex2} }
  \label{tab_18}
\begin{tabular}{ |p{1.7cm}||p{1.5cm}|p{1.5cm}| p{1.7cm}|| p{1.5cm}|p{1.5cm}|}
 \hline
    \multicolumn{3}{|c|}{Method via Hankel pencil} & \multicolumn{3}{|c|}{AAA method}\\
   \hline
$\sigma=10^{-7}$ & $e(\boldsymbol{z})$  & $e(\gamra) $ & $\sigma=10^{-7}$ & $e(\boldsymbol{z})$  &$e(\gamra) $  \\[2mm]
 \hline
min & $4.35$e--$03$    &$7.99$e--$05$ & min &$1.26$e--$01$ & $4.89$e--$03$  \\[2mm]
max  & $3.00$e--$02$    &$5.66$e--$04$ & max &$2.03$e+$00$ & $7.86$e--$02$   \\[2mm]
average & $1.42$e--$02$    &$2.66$e--$04$ & average & $1.48$e+$00$ & $4.44$e--$02$   \\[2mm]
 \hline
 \end{tabular}
  \qquad 
  \begin{tabular}{ |p{1.7cm}||p{1.5cm}|p{1.5cm}|p{1.7cm}|| p{1.5cm}|p{1.5cm}|}
 \hline
    \multicolumn{3}{|c|}{Method via Hankel pencil} & \multicolumn{3}{|c|}{AAA method}\\
   \hline
$\sigma=10^{-8}$ & $e(\boldsymbol{z})$ & $e(\gamra) $ & $\sigma=10^{-8}$ & $e(\boldsymbol{z})$  &$e(\gamra) $  \\[2mm]
 \hline
min & $7.89$e--$04$    &$4.48$e--$06$ & min & $2.6$e--$02$ & $1.01$e--$03$  \\[2mm]
max  & $1.28$e--$03$    &$5.58$e--$05$ & max & $2.6$e--$01$ & $1.01$e--$02$   \\[2mm]
average & $2.98$e--$03$    &$2.4$e--$05$ & average & $1.25$e--$01$ & $4.84$e--$03$   \\[2mm]
 \hline
 \end{tabular}
 \end{table}

\begin{example}\label{ex4}

   We consider again recovery of rational functions $r$ from Examples \ref{ex3}--\ref{ex2}. In Tables \ref{tab_16}--\ref{tab_18} we present the results for $N=15$, $L=20$ (for the modified method $L\leq 2N$) and different values of $\sigma$. 

  Analyzing results of Tables \ref{tab_16}--\ref{tab_18}, we conclude that for $r$ from Example \ref{ex1} both methods work very well, while for   $r$ from Examples \ref{ex3}  and \ref{ex2} we get a slightly better recovery  applying our method. This is because the Hankel pencil-based approach is more effective for the recovery of poles with certain locations.
 For Example \ref{ex2} the poles $z_1^{(-)}=0.2$, $z_2^{(-)}=0.5$ and $z_1^{(+)}=2$  are reconstructed equally good by both methods, while the Hankel pencil method shows better performance for the recovery of the pole  $z_2^{(+)}=50$ than the AAA method (see Table   \ref{tab_14}).

\begin{table}[h!]
\centering
 \caption{\small Reconstruction errors $|z_j-\tilde{z}_j|$ for  poles $z_j$ from Example \ref{ex2}   with the Hankel pencil and the AAA methods}
  \label{tab_14}
\begin{tabular}{ |p{1.8cm}||p{2.2cm}|p{2.2cm}| p{2.2cm}|p{2.2cm}|}
 \hline
    \multicolumn{5}{|c|}{Method via Hankel pencil} \\
   \hline
$\sigma=10^{-7}$ & $z_1^{(-)}=0.2$ & $z_2^{(-)}=0.5$  & $z_1^{(+)}=2$  &$z_2^{(+)}=50$  \\[2mm]
 \hline
min & $2.82$e--$08$    &$6.98$e--$08$ & $2.08$e--$07$ & $4.35$e--$03$  \\[2mm]
max  & $5.09$e--$07$    &$3.92$e--$07$ & $1.14$e--$06$ & $3$e--$02$   \\[2mm]
average & $2.51$e--$07$    &$1.79$e--$07$ & $2.08$e--$07$ & $1.72$e--$02$   \\[2mm]
 \hline
 \end{tabular}
 \begin{tabular}{ |p{1.8cm}||p{2.2cm}|p{2.2cm}| p{2.2cm}|p{2.2cm}|}
 \hline
\multicolumn{5}{|c|}{AAA method} \\
   \hline
$\sigma=10^{-7}$ & $z_1^{(-)}=0.2$ & $z_2^{(-)}=0.5$  & $z_1^{(+)}=2$  &$z_2^{(+)}=50$  \\[2mm]
 \hline
min & $1.15$e--$07$    &$6.42$e--$08$ & $6.65$e--$07$ & $1.26$e--$01$  \\[2mm]
max  & $7.4$e--$07$    &$5.24$e--$07$ & $2.55$e--$06$ & $2.03$e+$00$   \\[2mm]
average & $4.63$e--$07$    &$2.34$e--$07$ & $1.19$e--$06$ & $1.14$e+$00$   \\[2mm]
 \hline
 \end{tabular}
  \end{table}

\end{example}

\section*{Acknowledgement}
 The author would like to thank Ion Victor Gosea (MPI Magdeburg) for several valuable discussions regarding the  sensitivity analysis. The author would also like to thank the anonymous referees for their valuable comments that helped to improve the manuscript.

\end{document}